\documentclass[11pt]{article}

\usepackage[round]{natbib}

\usepackage{geometry}

\usepackage{graphicx}

\usepackage{amsmath}
\usepackage{amsfonts}
\usepackage{enumerate}
\usepackage{lscape}
\usepackage{multirow}

\usepackage{amsthm}
\newtheorem{definition}{Definition}
\newtheorem{lemma}{Lemma}
\newtheorem{theorem}[lemma]{Theorem}
\newtheorem{corollary}[lemma]{Corollary}

\usepackage{algorithm}
\usepackage{algpseudocode}
\newcommand{\AlgorithmFontSize}{\footnotesize}

\makeatletter
\newcommand{\firstcline}[1]{%
  \multicolumn1c{%
    \global\backup@length\ht\@arstrutbox
    \global\advance\backup@length\dp\@arstrutbox
    \global\advance\backup@length\arrayrulewidth
     \raise\extratabsurround\copy\@arstrutbox
    }\\[-\backup@length]\cline{#1}
}
\makeatother

\usepackage{pdflscape}

\usepackage{url}

\usepackage{pgfplots}
\usepackage{pgfplotstable}
\pgfplotsset{compat=1.11,
        /pgfplots/ybar legend/.style={
        /pgfplots/legend image code/.code={%
        \draw[##1,/tikz/.cd,bar width=3pt,yshift=-0.2em,bar shift=0pt]
                plot coordinates {(0cm,0.8em)};},
},
}
\newcommand{\change}[1]{{#1}}
\newcommand{\abstractchange}[1]{{#1}}

\usepackage{authblk}

\begin{document}

\title{When Lift-and-Project Cuts are Different}

\author[1]{Egon Balas} 
\author[2]{Thiago Serra}

\affil[1]{Carnegie Mellon University}
\affil[2]{Bucknell University}

\date{}

\maketitle

\begin{abstract} 
In this paper, 
we present a method to determine if a lift-and-project cut for a mixed-integer linear program is \abstractchange{irregular}, 
in which case the cut is \abstractchange{not} equivalent to \abstractchange{any} intersection cut from \abstractchange{the bases} of the linear relaxation. 
This is an important question  
due to the intense research activity for the past decade on cuts from multiple rows of simplex tableau 
as well as on lift-and-project cuts from non-split disjunctions. 
While it is known since \cite{BalasP} that 
\abstractchange{lift-and-project cuts from split disjunctions 
are always equivalent to intersection cuts and consequently to such multi-row cuts, }
\cite{BalasK} have recently shown that there is a necessary and sufficient condition in the case of arbitrary disjunctions:  
a lift-and-project cut is regular if, and only if, it corresponds to a regular \abstractchange{basic} solution of the Cut Generating Linear Program~(CGLP).
This paper has four contributions. First, we state a result that simplifies the verification of regularity for basic CGLP solutions from \cite{BalasK}. 
Second, we 
provide a 
mixed-integer formulation that checks whether there is a regular CGLP solution for a given cut \abstractchange{that is regular in a broader sense, which also encompasses irregular cuts that are implied by the regular cut closure}. Third, we describe a numerical procedure \abstractchange{based on such formulation} that \abstractchange{identifies irregular lift-and-project cuts}. Finally, 
we use this method to evaluate how often lift-and-project cuts \abstractchange{from simple $t$-branch split disjunctions} are irregular, 
and thus not equivalent to multi-row cuts, 
on 74 instances of the MIPLIB benchmarks.

\end{abstract}

\section{Introduction}

Many techniques to generate cutting planes for a Mixed-Integer Linear Program~(MILP) 
are equivalent to one another under certain conditions. 
Since some are more general and usually more expensive computationally, 
it is important to determine if and when they generate cuts that others cannot. 
In this paper, 
we introduce a technique to verify if a lift-and-project cut~\citep{CGLP,CGLP2} from an arbitrary disjunction does not 
correspond to a standard intersection cut~\citep{IC} 
and analyze computational results on several instances for insights. 
This study is particularly relevant due to the recent research activity around intersection cuts from multiple rows of the simplex tableau since their introduction by \cite{TwoRows}, 
which in the case of $q$ rows are equivalent to disjunctive cuts from a $2^q$-term disjunction~\citep{MRDP}. 
More specifically, 
it helps us find out about the converse: \change{how often lift-and-project cuts from multi-term disjunctions cannot be directly obtained as intersection cuts from solutions of a simplex tableau, and therefore are different. We test when these different cuts -- which are denoted as irregular -- are found in practice, what may be inducing their occurrence, and how strong they are.}

We can obtain a lift-and-project cut by solving a Cut Generating Linear Program~(CGLP), 
which defines valid inequalities for an MILP relaxation consisting of a disjunctive program~\citep{Balas79,DP}. 
Such \change{disjunctive programs} are often unions of disjoint polyhedra that exclude the region around a particular solution $x=\bar{x}$ of the Linear Program~(LP) relaxation.  
The most common \change{ disjunctive program} consist of intersecting the LP relaxation 
with a \emph{split disjunction} of the form~$\{ x : \pi x \leq \pi_0 \} \vee \{ x : \pi x \geq \pi_0 + 1 \}$ 
when $\pi_0 < \pi \bar{x} < \pi_0+1$ and no MILP solution is removed by the disjunction. 
For that case, 
\citet{BalasP} have shown that there is a correspondence between lift-and-project cuts from basic CGLP solutions 
and intersection cuts from basic solutions of the LP relaxation, feasible or not. 
That entails a more efficient procedure to implicitly solve CGLPs from split disjunctions by pivoting among LP bases, 
which has been implemented in a number \change{ of} solvers including CglLandP~\citep{BalasB} in COIN-OR\footnote{Available at \url{http://www.coin-or.org}}.

The equivalence identified by \citet{BalasP} 
associates  
lift-and-project cuts from a \emph{simple} split disjunction of the form~$\{ x : x_k \leq 0 \} \vee \{ x : x_k \geq 1 \}$ 
with Gomory Mixed-Integer cuts~\citep{GMI} from the row of some  simplex tableau, feasible or not, defining the value of $x_k$ in terms of nonbasic continuous variables.  
Similarly, 
strengthening those lift-and-project cuts by 
modularizing 
the coefficients associated with integer nonbasic variables~\citep{Strengthening,CGLP} 
makes them equivalent to GMI cuts from the corresponding row of some simplex tableau that defines $x_k$ in terms of nonbasic variables, 
some or all of which may be integer-constrained.  

More recently, 
\citet{BalasK} have shown that the correspondence between lift-and-project cuts and intersection cuts 
may not necessarily hold for 
lift-and-project cuts from arbitrary disjunctions. 
More specifically, 
they have proven that it holds if, and only if, 
there is a basic CGLP solution associated with the cut 
that maps to an LP basis that corresponds to a standard intersection cut. 
These basic CGLP solutions mapping to LP~bases are called \emph{regular}. 
A lift-and-project cut is then called \emph{regular} if there exists a corresponding regular basic CGLP solution   
and \emph{irregular} otherwise. 

The elegance and convenience of generating intersection cuts from the simplex tableau 
has motivated a recent stream of theoretical work on generating cuts from two rows of the simplex tableau~\citep{TwoRows,CornMargot} and subsequently more rows~\citep{LatticeBorozan,LatticeBasu}, 
on how these cuts can be strengthened when nonbasic variables are integer~\citep{DeyLifting,ConfortiLifting,BasuLifting,FukasawaLifting}, 
and several other variants. 
The reader is referred to \citet{CornerReview} and \citet{GeometryReview} for a broader review of this line of work, 
which has been accompanied by extensive computational experimentation~\citep{MultiRowsExp,BasuExp,DeyExp,Louveaux}.

However, there are other ways in which one can exploit that more than one integer variable is fractional in $\bar{x}$. 
For example, we can generate lift-and-project cuts using disjunctions with more than two terms, 
and those may yield irregular cuts instead. 
A natural generalization of the commonly used split disjunction is defined by~\cite{LiRichard} as the $t$-\emph{branch split disjunction} 
\[
\bigvee\limits_{S \subseteq \{1,2,\ldots,t\}} \big\{ x : \pi^i x \leq \pi^i_0, \text{ if } i \in S; ~ \pi^i x \geq \pi^i_0 + 1, \text{ if } i \notin S \big\}, 
\]
of which the cross disjunction~\citep{Cross} corresponds to the special case where $t=2$. 
\cite{CrossExps} observed that the resulting cuts close a substantially larger gap in comparison to split cuts. 
Furthermore, 
\cite{AndersenCL} as well as \cite{Kis} have shown examples of disjunctive cuts that do not correspond to intersection cuts.

\begin{figure}
    \centering
    \includegraphics[width=0.23\textwidth]{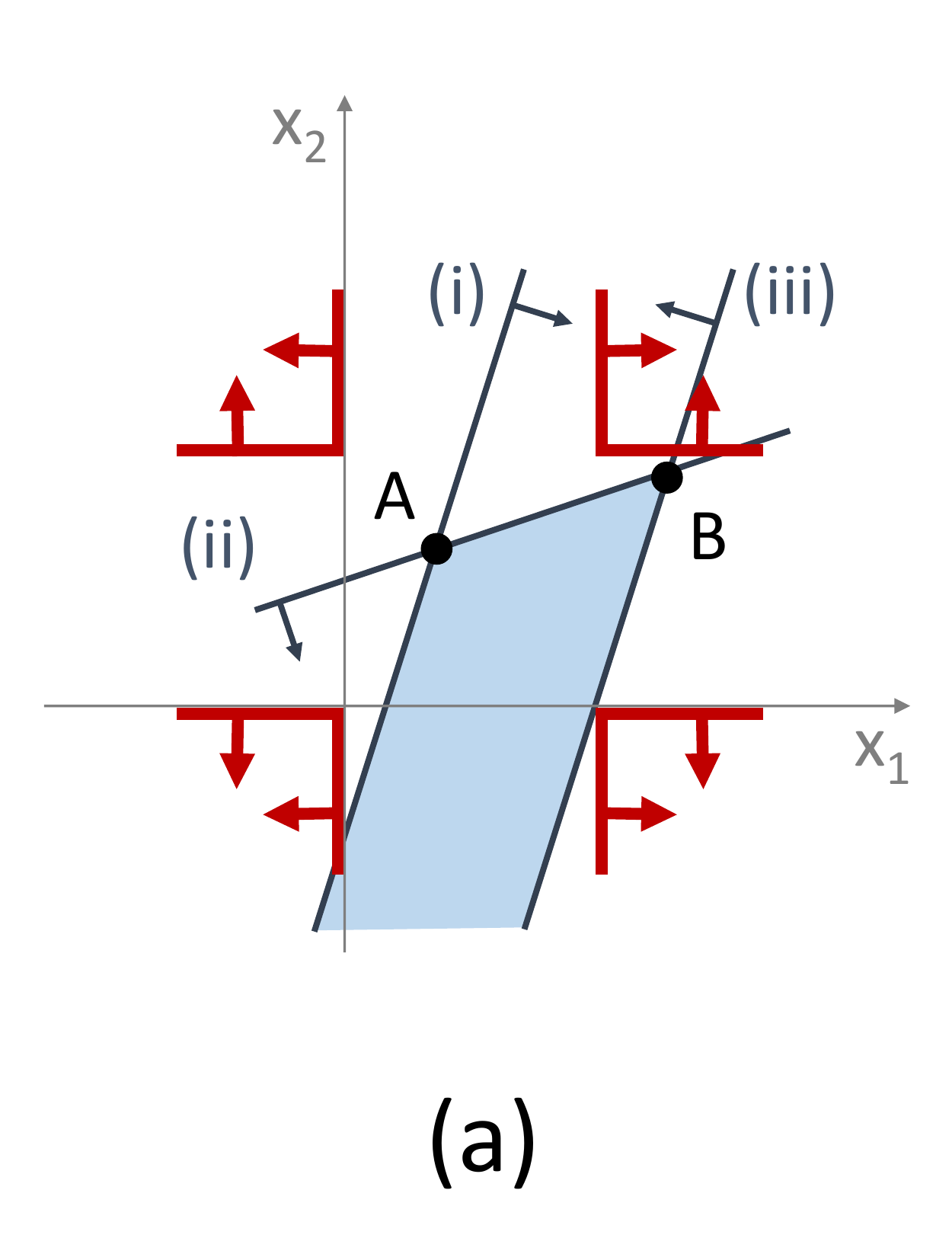}~~~\includegraphics[width=0.23\textwidth]{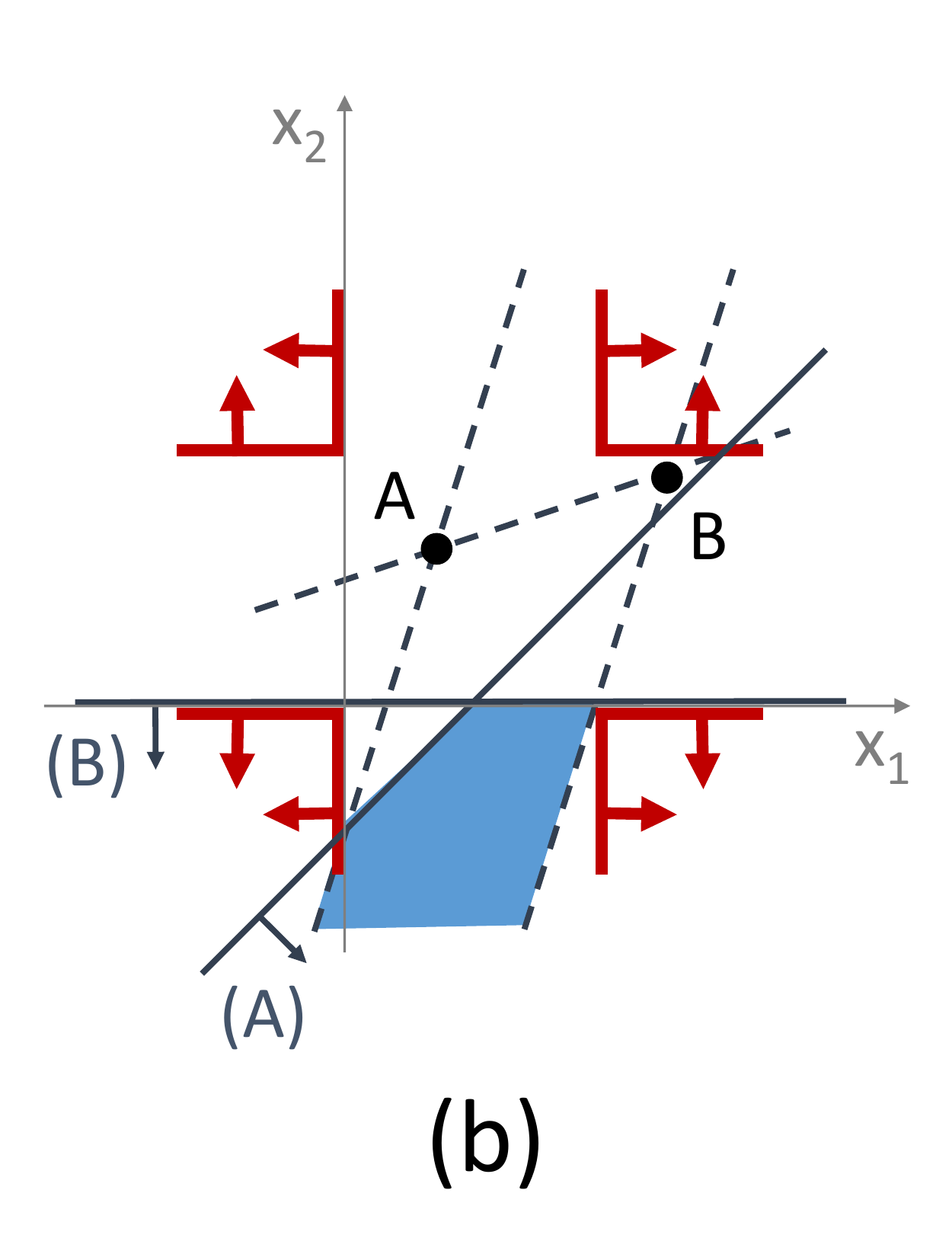}~~~\includegraphics[width=0.23\textwidth]{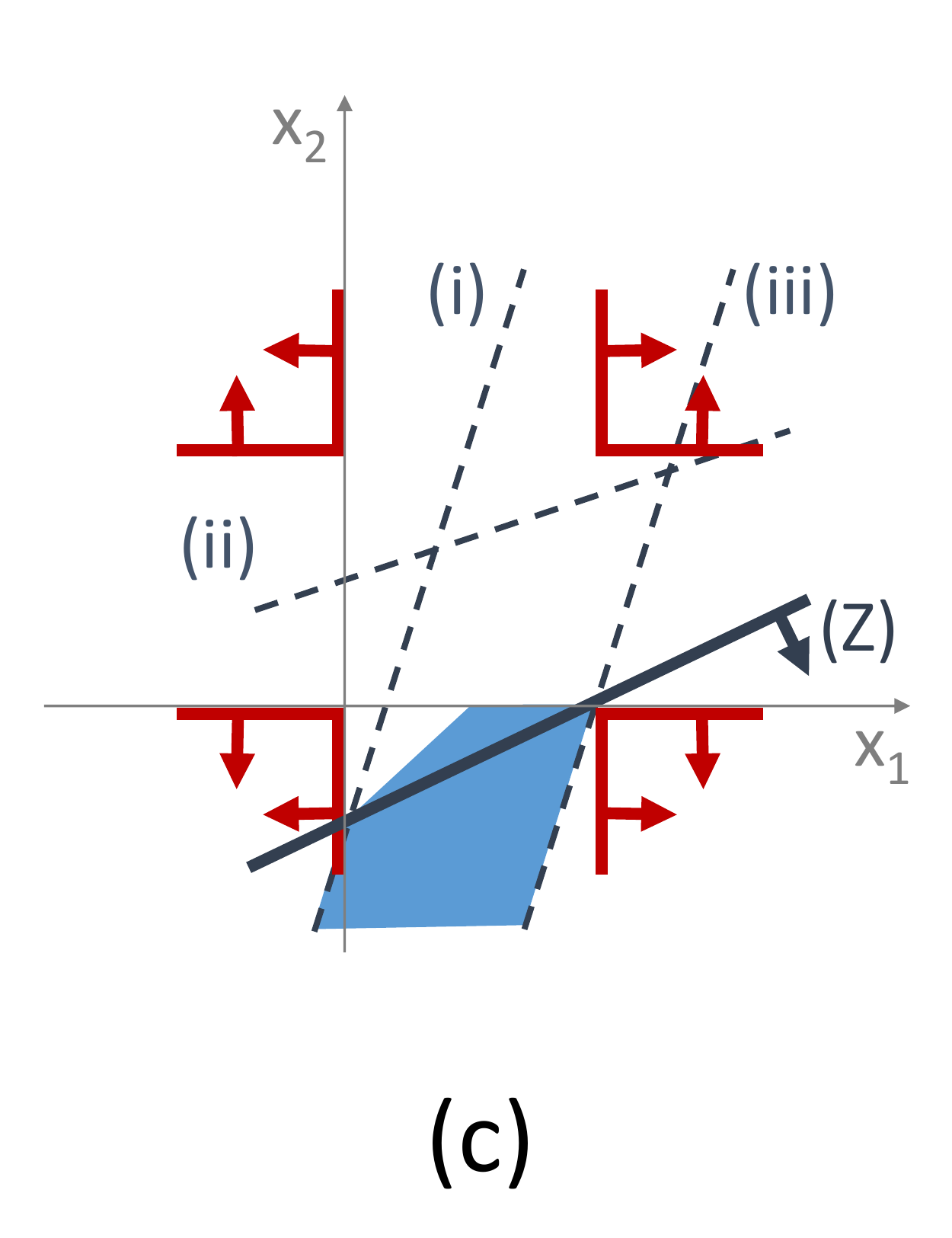}~~~\includegraphics[width=0.23\textwidth]{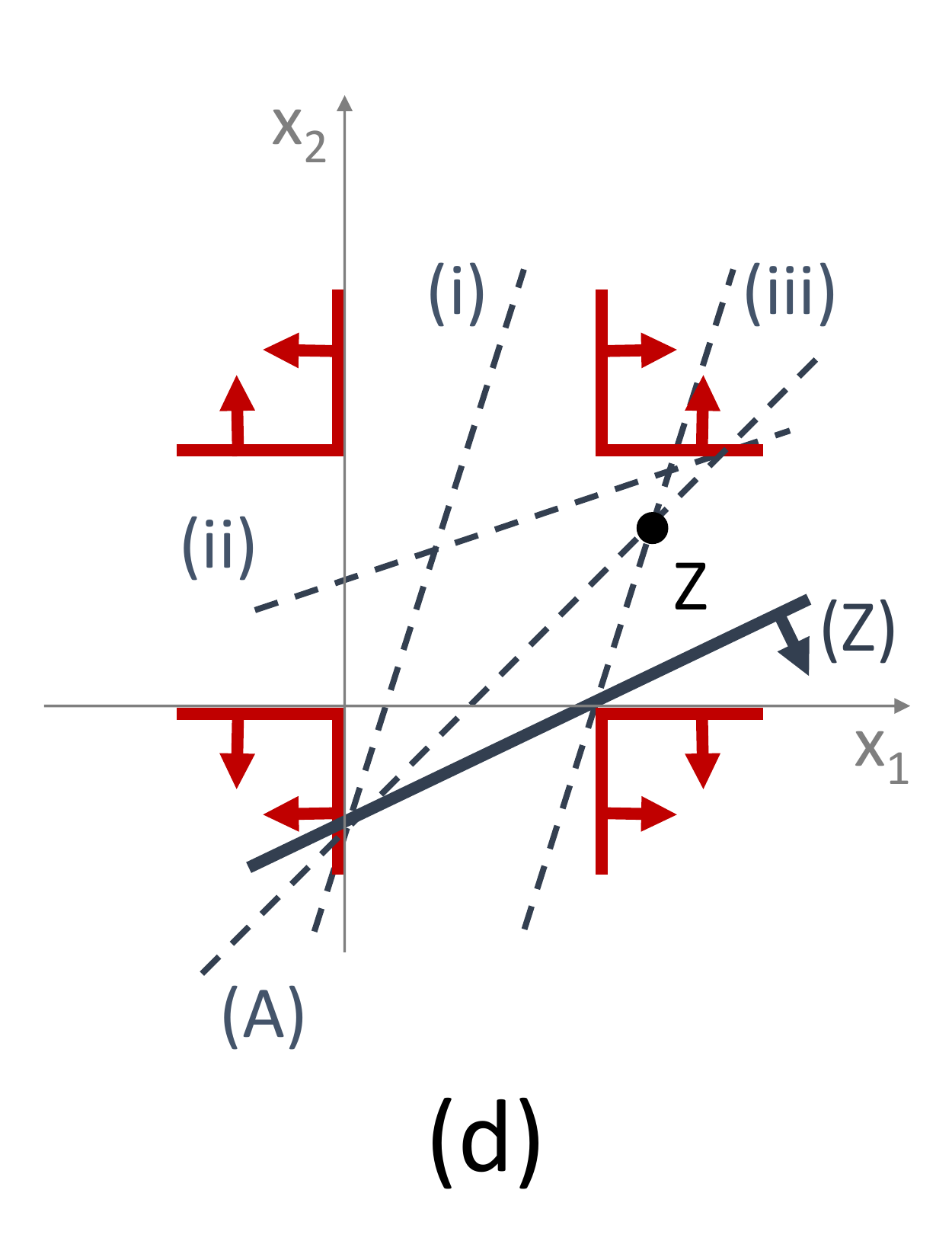}
    \caption{(a) LP relaxation and a 2-branch split disjunction; (b) Non-dominated intersection cuts; 
    (c) Irregular cut dominating all intersection cuts; and (d) Additional inequality making previous cut regular.}
    \label{fig:illustration}
\end{figure}

Figure~\ref{fig:illustration} illustrates another such case. In Figure~\ref{fig:illustration}a, we have an LP relaxation defined by inequalities (i) $6 x_1 - 2 x_2 \geq 1$, (ii) $2 x_1 - 6 x_2 \geq -3$, and (iii) $-3 x_1 + x_2 \geq -3$ in lighter blue and a 2-branch split disjunction $\{x : x_1 \geq 1, x_2 \geq 1 \} \vee \{ x: x_1 \leq 0, x_2 \geq 1\} \vee \{ x: x_1 \leq 0, x_2 \leq 0\} \vee \{ x: x_1 \geq 1, x_2 \leq 0 \}$ in red. The LP has basic feasible solutions at points $A = (\frac{3}{8},\frac{5}{8})$ and $B = (\frac{21}{16},\frac{15}{16})$. In Figure~\ref{fig:illustration}b, we have the non-dominated intersection cuts (A) $2 x_1 - 2 x_2 \geq 1$ and (B) $x_2 \leq 0$ defining the regular cut closure in the darker blue area. These cuts are respectively associated with the LP bases in points A and B. In Figure~\ref{fig:illustration}c, we show another valid inequality (Z) $x_1 - 2 x_2 \geq 1$ 
cuts through the regular cut closure and in fact dominates the intersection cuts from the LP bases within the LP relaxation. While each intersection cut can be obtained with the two inequalities that are satisfied at equality of the corresponding LP basis, the three inequalities are necessary in the case of (Z): we need both inequalities (ii) and (iii) to certify that (Z) does not separate any point at the intersection of the LP relaxation with the first term of the disjunction, 
and similarly we need inequality (i) in the case of the the third term of the disjunciton.
In Figure~\ref{fig:illustration}d, we show that (Z) would be an intersection cut from an LP basis at $Z = (\frac{5}{4},\frac{3}{4})$ if the LP relaxation also includes the valid inequality (A). Hence, we cannot infer that a lift-and-project cut is irregular by only inspecting the CGLP solution from which we obtained the cut, since it could also be generated in other ways.

In this work, we focus on the equivalence with respect to lift-and-project cuts without strengthening, 
%
and we present 
the following contributions. 
First, 
we state a result that simplifies the verification of regularity for basic CGLP solutions from~\citet{BalasK} 
and show 
how it extends 
to CGLP solutions that are not basic  
in Section~\ref{sec:cglp_reg}. 
Second, 
we introduce and prove the validity of an MILP that checks whether there is a regular CGLP solution, \change{basic or not}, for a given cut in Section~\ref{sec:cut_reg}. 
Third, 
we describe a numerical procedure based on such MILP that verifies if a lift-and-project cut is irregular \change{in a stricter sense} in Section~\ref{sec:algo}.
Finally, 
we present computational results from 74 benchmark instances in Section~\ref{sec:exps}, 
and we use these results to analyze what factors may lead to a higher incidence of \change{strictly}  irregular cuts and how these cuts compare with the other cuts in Section~\ref{sec:disc}. 

\section{Preliminaries}\label{sec:prelim}

Let us consider a mixed 0$-$1 linear program with rational coefficients 
\begin{align*}\tag{$\mathcal{P}$}\label{eq:mip}
\min \Big\{ cx ~ : ~ Ax \geq b, x\geq 0, ~ x_j \in \{0, 1\}, ~ j=1, \ldots, p \Big\}, 
\end{align*}
where $A$ is an $m \times n$ matrix. 
Let $\tilde{A} x \geq \tilde{b}$ denote the constraint set of the LP relaxation $P := \{ x : Ax \geq b, x \geq 0, x_j \leq 1, j=1, \ldots, p\}$ 
and let $P_I := \{ x : x \in P, x_j \in \{0,1\}, j = 1, \ldots, p\}$ denote the feasible set of~\eqref{eq:mip}.  
Hence, $\tilde{A}$ is a $q \times n$ matrix, 
where $q = m + n + p$.
Let $Q := \{1, \ldots, q\}$.

Furthermore, 
let $\bar{x}$ be a basic optimal solution of the LP relaxation of~\eqref{eq:mip}, 
i.e., $\bar{x} = \arg \min \{c x : \tilde{A} x \geq \tilde{b} \}$. 
If $\bar{x}$ is not feasible for~\eqref{eq:mip}, 
we can define a disjunction \change{ $\vee_{t \in T} D^t x \geq d_0^t$} 
that contains the feasible set of~\eqref{eq:mip} and not $\bar{x}$. 
For example, if there is a nonempty set $K \subseteq \{1, \ldots, p\}$ for which $0 < \bar{x}_k <1$ for every $k \in K$, 
then 
we can define a disjunction of the form 
\begin{equation}\label{eq:disj_k}
\bigvee_{K' \subseteq K} \left(
\begin{array}{ll}
x_k \geq 1, & k \in K' \\
x_k \leq 0, & k \in K \setminus K'
\end{array}
\right), 
\end{equation}
which 
we denote 
as a \emph{simple} $t$-branch split disjunction\footnote{Unless noted otherwise, we will use $t$ to index terms of a  disjunctive program in general form instead of parameterizing a $t$-branch split disjunction. We will use a set $K$ to denote a particular simple $t$-branch split disjunction.}, where $t = |K|$.  
More generally, we can generate lift-and-project cuts by intersecting sets such as~\eqref{eq:disj_k} with the linear relaxation of~\eqref{eq:mip}. 
More generally, we have a \change{ disjunctive program} of the form 
\begin{equation}\label{eq:disj}
\bigvee_{t \in T} \left(
\begin{array}{c}
\tilde{A} x \geq \tilde{b} \\
D^t x \geq d_0^t
\end{array}
\right),
\end{equation}
\change{ where, without loss of generality, $D^t$ is a $r \times n$ matrix. }
The classic formulation for the Cut Generating Linear Program~(CGLP) used to find a lift-and-project cut \change{ among valid inequalities for}~\eqref{eq:disj} separating $\bar{x}$, which excludes some dominated inequalities, is as follows: 
\change{
\begin{alignat}{8}
\min ~ & \alpha \bar{x} & - & \beta \label{cglp:obj} \\
& \alpha &&& -u^t \tilde{A} && - v^t D^t & = & ~ 0, & ~ ~ ~ t \in T \label{cglp:alpha} \\
&&& \beta & -u^t \tilde{b} && - v^t d^T_0 & = & ~ 0, & ~ ~ ~  t \in T \label{cglp:beta} \\
&&& & \sum_{t \in T} \sum_{i=1}^q u^t_i && + \sum_{t \in T} \sum_{i=1}^r v^t_i  & ~= & ~ 1 \label{cglp:norm} \\
&&& & u^t, && v^t & \geq & ~ 0, & ~ ~ ~  t \in T \label{cglp:nonneg} \\
& \alpha \in \mathbb{R}^n, && \beta \in \mathbb{R} \label{cglp:urs}
\end{alignat} 
}
\change{In the formulation above, the data consist of $\tilde{A}$, $\tilde{b}$, $D^t$ and $d^t_0$ for $t \in T$. The variables are $\alpha$, $\beta$, $u^t$ and $v^t$ for $t \in T$. If $\alpha \bar{x} - \beta < 0$, then $\alpha x \geq \beta$ is a valid inequality separating $\bar{x}$.}

\change{For every choice of multipliers $\{ v^t \}_{t \in T}$ such that $v^t \geq 0$ and $v^t \neq 0$, where $v^t$ is a vector in which each element corresponds to an inequality of $D^t x \geq d_0^t$, we can define a polyhedron $S\left( \{ v^t \}_{t \in T} \right) := \{ x : (v^t D^t) x \leq v^t d_0^t \}$. 
If $(v^t D^t) \bar{x} < v^t d_0^t$ for each $t \in T$, 
then such polyhedron contains $\bar{x}$ as an interior point but no point in $P_I$, and thus it is possible to use polyhedron $S\left( \{ v^t \}_{t \in T} \right)$ to derive an intersection cut separating $\bar{x}$ from $P_I$ as discussed next.}

If we assume, without loss of generality, that the upper bounds on the binary variables are contained in $Ax \geq b$, 
then we extend $x$ with surplus variables of the form $x_{n+i} = \sum_{j =1}^n a_{i j} x_j - b_i$ for $i = 1, \ldots, m$  
and define the following LP cone $C(J)$ from each basic solution $x(J)$ of $P$:
\begin{align}
x_i  =  \bar{a}_{i 0} - \sum_{j \in J} \bar{a}_{i j} x_j, & \qquad i \in I \\
x_i  \geq  0, & \qquad i \in J 
\end{align}
where $I$ is the index set of basic variables and $J$ of the nonbasic variables. 
Cone $C(J)$ has an extreme ray $r^j(J)$ corresponding to each nonbasic variable $x_j, j \in J$, 
with $r^j_j(J) = 1$, $r^j_i(J) = - \bar{a}_{i j}$ for $i \in I$, and $r^j_i(J) = 0$ for $i \in J \setminus \{ j \}$. 
If there is a convex set $S$ containing $x(J)$ but no feasible integer point in its interior, 
which we denote as $P_I$-free, 
we can define the intersection cut~\citep{IC}  
\begin{align}
\sum_{j \in J} \frac{1}{\lambda^*_j} x_j \geq 1
\end{align}
separating $x(J)$ from $P_I$, 
where $\lambda^*_j$ defines the point at which each ray $x(J) + \lambda_j r_j(J)$ intersects the boundary of $S$. 
If some ray never intersects the boundary, 
then the corresponding coefficient of the intersection cut is zero instead of $\lambda^*_j$ inversed. 

\change{ If a lift-and-project cut $\alpha x \geq \beta$ separating $\bar{x}$ is obtained from a basic CGLP solution defined on a split disjunction $\{x : \pi x \leq \pi_0 \} \vee \{x: \pi x \geq \pi_0 + 1\}$, 
then the rows of the LP relaxation associated with the basic multipliers in the CGLP solution 
define the LP basis from which $\alpha x \geq \beta$ can be obtained as an intersection cut using the $P_I$-free set $\{x : \pi_0 \leq \pi x \leq \pi_0 +1 \}$. In this particular case, \citet{BalasP} have shown that there are exactly $n$ basic multipliers in $u$ associated with distinct rows of the LP relaxation.}

\change{For arbitrary disjunctions, the lift-and-project cuts that are equivalent to intersection cuts are always associated with at least one basic CGLP solution that generalizes the structure found in the case of split disjunctions. The following results, which are adapted from \citet{BalasK}, present the necessary and sufficient conditions for such equivalence:}

\begin{theorem}[\citet{BalasK}, Th. 10]\label{thm:suff}
\change{Let $(\bar{\alpha}, \bar{\beta}, \{ \bar{u}^t, \bar{v}^t \}_{t\ \in T})$ be a basic feasible solution to \eqref{cglp:alpha}--\eqref{cglp:urs} such that $\sum\limits_{i=1}^r \bar{v}^t_i>0$ for all $t \in T$. If there exists a nonsingular $n \times n$ submatrix $\tilde{A}_J$ of $\tilde{A}$ such that $\bar{u}^t_j = 0$ for all $j \notin J$ and $t \in T$, then the lift-and-project cut $\bar{\alpha} x \geq \bar{\beta}$ is equivalent to the intersection cut $\pi x_J \geq 1$ from the set \[
S(\bar{v}) := \{ x \in \mathbb{R}^n : (\bar{v}^t D^t) x \leq \bar{v}^t d^t_0, t \in T\}
\]
and the LP simplex tableau with nonbasic set $J$.}
\end{theorem}

\begin{theorem}[\citet{BalasK}, Th. 12]
\change{Let $\bar{w} = (\bar{\alpha}, \bar{\beta}, \{ \bar{u}^t, \bar{v}^t \}_{t\ \in T})$ be a basic feasible solution to \eqref{cglp:alpha}--\eqref{cglp:urs} such that $\sum\limits_{i=1}^r \bar{v}^t_i>0, t \in T$.}

\change{If $\bar{w}$ does not satisfy the (sufficient) condition of Theorem~\ref{thm:suff}, and there is no basic feasible solution $\tilde{w}$ to \eqref{cglp:alpha}--\eqref{cglp:urs} with $(\tilde{\alpha}, \tilde{\beta}) = \mu (\bar{\alpha}, \bar{\beta})$ for some $\mu > 0$ that satisfies the condition of Theorem~\ref{thm:suff}, then there exists no intersection cut from any member of the family of polyhedra
\[
S(v) := \{ x \in \mathbb{R}^n : (v^t D^t) x \leq v^t d^t_0, t \in T\}
\]
where $v \geq 0, v \neq 0$, equivalent to $\bar{\alpha} \bar{x} \geq \bar{\beta}$. }
\end{theorem}

\change{In other words, a cut $\bar{\alpha} x \geq \bar{\beta}$ from a regular basic solution $(\bar{\alpha}, \bar{\beta}, \{ \bar{u}^t, \bar{v}^t \}_{t\ \in T})$ is equivalent to the intersection cut from the LP cone associated with the cobasis indexed by $J$ and the $P_I$-free convex set defined by $\{ x : (\bar{v}^t D^t) x \leq \bar{v}^t d^t_0, t \in T \}$. 
More specifically, 
a positive CGLP multiplier $\bar{u}_i^t$ for some row $i$ of $Ax \geq b$ for any $t \in T$ maps the corresponding surplus variable $s_i$ as nonbasic in the LP, 
a positive multiplier $\bar{u}_{m+j}^t$ for the bound  $x_j \geq 0$ maps $x_j$ as nonbasic at the lower bound in the LP, 
and a positive CGLP  $\bar{u}_{m+p+j}^t$ for the bound $x_j \leq 1$ maps $x_j$ as nonbasic at the upper bound in the LP. 
The 
inequality $\bar{\alpha} x \geq \bar{\beta}$ only cuts off part of $P$ if $\sum\limits_{i=1}^r\bar{v}^t_i > 0$ for each $t \in T$, 
since otherwise $\bar{\alpha} x \geq \bar{\beta}$ is implied by $\tilde{A} x \geq \tilde{b}$, 
as shown in 
Lemma~1 of~\citet{BalasP}.
Furthermore, the sufficient condition for the regularity of $\bar{w}$ is also necessary for $\bar{\alpha} x \geq \bar{\beta}$, 
i.e., if the condition is not met by any other CGLP basis, then the cut is irregular. }

\change{The following definition for what constitutes a regular basic CGLP solution and regular lift-and-project cut 
is also adapted from~\citet{BalasK}:}

\begin{definition}[Regularity of bases and cuts~\citep{BalasK}]
\change{A feasible basis for the CGLP system \eqref{cglp:alpha}--\eqref{cglp:urs} and the associated solution will be called \emph{regular} if 
the basis satisfies the condition of Theorem~\ref{thm:suff}, and \emph{irregular} otherwise.
%
A cut defined by an irregular solution $w$ is \emph{irregular}, 
unless there exists a regular solution $w'$ with the same $(\alpha, \beta)$-component -- upon scaling by a positive multiplier -- as that of $w$, in which case the cut is \emph{regular}.}
\end{definition}

\section{Regularity of CGLP Solutions}\label{sec:cglp_reg}



The next Theorem gives a simple criterion for deciding whether a basic CGLP solution $\bar{w}$ is regular or not.

\begin{theorem}\label{thm:positive}
For a basic CGLP solution $\bar{w} = (\bar{\alpha}, \bar{\beta}, \{ \bar{u}^t, \bar{v}^t \}_{t \in T})$, 
let $\tilde{A}_N$ be the $|N| \times n$ submatrix of $\tilde{A}$ whose rows are indexed by 
$N(\bar{u}) := \{ j \in Q : \bar{u}_j^t > 0 \text{ for some } t \in T \}$. 
Then $\bar{w}$ is a regular solution if, and only if, 
$\tilde{A}_N$ is of full row rank.
\end{theorem}
\proof{Proof.}
\textit{Sufficiency}.  
Assume $\text{rank}(\tilde{A}_N) = |N|$. Then $|N| \leq n$. 
We show that in this situation $\bar{w}$ is regular.

\textit{Case 1}: $|N| = n$. Then $\tilde{A}_N$ is an $n \times n$ nonsigular submatrix of $\tilde{A}$ such that $u_j^t = 0$ for all $j \notin N$ and all $t \in T$, i.e., $\bar{w}$ satisfies the condition of \change{Theorem~\ref{thm:suff}}.

\textit{Case 2}: $|N| < n$. Then $\tilde{A}$ has $n - |N|$ rows $\tilde{A}_j$ with $u_j^t = 0$ which can be added to $\tilde{A}_N$ in order to form an $n \times n$ nonsingular matrix $\tilde{A}_{N'}$ since $\tilde{A}$ contains $I_n$. 
Substituting $\tilde{A}_{N'}$ for $\tilde{A}_N$ then reduces this case to Case 1. 

\textit{Necessity}. 
Assume $\text{rank}(\tilde{A}_N) < |N|$. We show that in this case $\bar{w}$ is irregular. 
In particular, 
any $n \times n$ nonsingular submatrix $\tilde{A}_J$ of $\tilde{A}$ has among its rows at most $\text{rank}(\tilde{A}_N)$ rows of $\tilde{A}_N$, 
thus leaving $|N| - \text{rank}(\tilde{A}_N)$ rows $j$ such that $\bar{u}_j^t > 0$ for some $t \in T$ outside of $\tilde{A}_J$. 
Therefore no such $\tilde{A}_J$ meets the condition of 
\change{Theorem~\ref{thm:suff}},  
hence $\bar{w}$ is irregular.
\endproof

In Figure~\ref{fig:illustration}c, note that the irregular cut (Z) can only be obtained by using the parallel inequalities (i) and (iii), for which reason the submatrix is not of full row rank.

\begin{corollary}\label{cor:all_regular}
If a CGLP solution $\bar{w} = (\bar{\alpha}, \bar{\beta}, \left\{ \bar{u}^t, \bar{v}^t \right\}_{t \in T} )$ is not basic but satisfies the conditions of Theorem~\ref{thm:positive} for regularity,  
then $\bar{\alpha} x \geq \bar{\beta}$ is valid for the closure of regular cuts. 
\end{corollary}
\proof{Proof.}
If $\bar{w}$ is not basic, then we can describe it as a proper convex combination of a set of basic CGLP solutions. 
Let these solutions be indexed by a given set $\mathcal{B}$, 
so that 
$(\bar{\alpha}, \bar{\beta}, \left\{ \bar{u}^t, \bar{v}^t \right\}_{t \in T} ) =$ $\sum_{b \in \mathcal{B}} \lambda_b \left( \tilde{\alpha}^b, \tilde{\beta}^b, \left\{ \tilde{u}^{t^b}, \bar{v}^{t^b} \right\}_{t \in T} \right)$, $\lambda > 0$, and $\sum_{b \in \mathcal{B}} \lambda_b = 1$.

Note that $\tilde{u}^{t^b}_j > 0$ implies that $\bar{u}^t_j > 0$, 
and thus $N(\tilde{u}^b) \subseteq N(\bar{u})$ for all $b \in \mathcal{B}$. 
Since the submatrix of $\tilde{A}$ on the rows of $N(\bar{u})$ is of full row rank, 
then there exists a set $N$ such that $|N| = n$ and $N(\bar{u}) \subseteq N$ for which $\text{rank}(\tilde{A}_N) = n$, 
and thus $N(\tilde{u}^b) \subseteq N$ for all $b \in \mathcal{B}$. 
That implies that all inequalities of the form $\tilde{\alpha}^b x \geq \tilde{\beta}^b$ for each $b \in\mathcal{B}$ 
define intersection cuts from a same basis, 
and thus each of those is regular. 
\endproof

\change{Motivated by Corollary~\ref{cor:all_regular}, we consider for the rest of the paper a broader definition of regularity for CGLP solutions, 
which conversely implies a stricter definition for cut irregularity 
that disregards certain cuts that are implied by the set of regular cuts:}

\begin{definition}[Extended Regularity]\label{def:extended}
\change{A feasible solution for the CGLP system \eqref{cglp:alpha}--\eqref{cglp:urs} will be called \emph{extended regular} if it satisfies all conditions of Theorem~\ref{thm:positive} other than being basic. A lift-and-project cut $\alpha x \geq \beta$ is \emph{strictly irregular} if 
there is no extended regular solution with the same $(\alpha, \beta)$-component upon scaling by a positive multiplier.}
\end{definition}

Thus deciding whether a basic CGLP solution $\bar{w}$ is regular or not is straightforward. 
\change{
Furthermore, 
any regular CGLP solution suffices to prove that a given cut is not strictly irregular. 
Next we examine how to use that to find strictly irregular cuts.}

\section{Regularity of Cuts from CGLP Solutions}\label{sec:cut_reg}

Given an irregular CGLP solution  $\bar{w} = (\bar{\alpha}, \bar{\beta}, \{ \bar{u}^t, \bar{v}^t \}_{t \in T})$, 
we define a mixed-integer program based on the CGLP to establish whether 
there is \change{an extended regular} CGLP solution 
$\tilde{w} = ( \tilde{\alpha}, \tilde{\beta}, \{ \tilde{u}^t, \tilde{v}^t \}_{t \in T} )$    
such that $(\tilde{\alpha}, \tilde{\beta}) = \theta (\bar{\alpha}, \bar{\beta})$ for some $\theta > 0$. 
In comparison to the CGLP formulation, 
we add a variable $\theta \in [0,1]$ 
and restrict the value of $(\alpha, \beta)$ to $\theta (\bar{\alpha}, \bar{\beta})$. 
Furthermore, we remove the normalization constraint~\eqref{cglp:norm} and introduce a binary upper bounding variable $\delta_j$ for each $u_j$ in order to model the set $N$ of indices $j \in Q$ such that $u_j^t > 0$ for some $t \in T$. 
We embed Theorem~\ref{thm:positive} by restricting the size of such set to at most $n$ 
and requiring that the submatrix $\tilde{A}_N$ to be of full row rank. 
We denote the resulting formulation as the Regular Cut Verifier MILP~(RCV-MILP):

\change{ 
\begin{alignat}{7}
\max ~ ~ ~ \theta ~ ~ \\
\theta \bar{\alpha} & & - u^t \tilde{A} & - v^t D^t & = & ~0, & t \in T \label{mip:alpha} \\
\theta \bar{\beta} & & - u^t \tilde{b} & - v^t d^t_0 & = & ~0, & t \in T \label{mip:beta} \\
& \delta_j & - u^t_j & & \geq & ~0, & j \in Q, ~ t \in T \label{mip:bin_limit} \\
& \sum_{j \in Q} \delta_j & & & \leq & ~n \label{mip:number} \\
& \sum_{j \in N} \delta_j & & & \leq & ~\text{rank}(\tilde{A}_N), & N \subseteq Q \label{mip:rank} \\
& & u^t, & v^t & \geq & ~0, & t \in T \label{mip:nonneg} \\
&&& \delta_j & \in & ~\{0, 1\}, & j \in Q \label{mip:bin_domain} \\ 
&&& \theta & \in & ~[0, 1] 
\end{alignat}
}
\change{ The data consist of $\bar{\alpha}$, $\bar{\beta}$, $\tilde{A}$, the rank of all submatrices of $\tilde{A}$, $\tilde{b}$, $D^t$ and $d^t_0$ for $t \in T$. The variables are $\theta$, $\delta$, $u^t$ and $v^t$ for $t \in T$. If $\theta > 0$, then $\bar{\alpha} x \geq \bar{\beta}$ is a regular cut.}

Some comments regarding RCV-MILP are in order. 
First, 
constraints  \eqref{mip:bin_limit} and \eqref{mip:bin_domain} 
define an implicit normalization constraint $\| u \|_{\infty} \leq 1$, 
which may prevent us from finding the cut in the same scale. 
Hence, variable $\theta$ is necessary even though normalization~\eqref{cglp:norm} is removed in comparison to the CGLP. 
\change{Since it is not immediate that the implicit normalization guarantees that RCV-MILP is always bounded, 
we set an upper limit of 1 to $\theta$.} 
Second, 
the equivalence between a RCV-MILP solution $(\breve{\theta}, \breve{\delta}, \{ \breve{u}^t, \breve{v}^t \}_{t \in T} )$ and a regular CGLP solution denoting a cut equivalent to $\bar{\alpha} x \geq \bar{\beta}$ is not direct. 
Instead of simply stating a corresponding CGLP solution $(\breve{\theta} \bar{\alpha}, \breve{\theta} \bar{\beta}, \{ \breve{u}^t, \breve{v}^t \}_{t \in T})$, 
we may need to first scale the RCV-MILP solution if \change{$\sum_{t \in T} \sum_{i=1}^r \breve{u}^t_i + \sum_{t \in T} \sum_{i=1}^r \breve{v}^t_i \neq 1$} in order to satisfy normalization~\eqref{cglp:norm}. 
Such scaling is done a number of times in the next proof. 
Last, 
\change{ 
for each subset $N \subseteq Q$ we assume that the rank of $\tilde{A}_N$ is given as an input.  
However, 
}
the number of subsets of $Q$, and consequently the number of rows due to constraint \eqref{mip:rank}, can be very large. 
In Section~\ref{sec:algo}, we address that by \change{ iteratively adding to the formulation only the relevant subsets of $Q$ and computing the rank of the corresponding submatrices of $\tilde{A}$}. 

The next result proves the validity of RCV-MILP. 
In what follows, 
we keep denoting the set of rows with positive multipliers as $N(\bar{u})$. 
Furthermore, 
let $\delta(\bar{u})$ be a vector in which $\delta_i = 1$ if $\bar{u}_i^t > 0 \text{ for some } t \in T$ and $\delta_i = 0$ otherwise.

\begin{theorem}
Let $\bar{w} = (\bar{\alpha}, \bar{\beta}, \{ \bar{u}^, \bar{v}^t \}_{t \in T})$ be a basic optimal solution of CGLP 
and let $(\check{\theta}, \check{\delta}, \{ \check{u}^t, \check{v}^t \}_{t \in T})$ be an optimal solution of RCV-MILP 
for cut $\bar{\alpha} x \geq \bar{\beta}$. 
Then the cut is \change{not strictly irregular} if, and only if, $\check{\theta} > 0$. 
\end{theorem}
\proof{Proof.}
If $\bar{w}$ is a regular CGLP solution, 
then $(1, \delta(\bar{u}), \{ \bar{u}^t, \bar{v}^t \}_{t \in T})$ is an optimal solution of RCV-MILP. 
First, 
$\theta = 1$ implies that 
constraints 
\eqref{cglp:alpha} and \eqref{cglp:beta} 
are equivalent to 
\eqref{mip:alpha} and \eqref{mip:beta}, 
whereas constraints \eqref{cglp:nonneg} and \eqref{mip:nonneg} are the same.  
Second,
normalization \eqref{cglp:norm} implies that $\bar{u} \leq 1$ 
and thus $\delta(\bar{u})$ satisfies constraints \eqref{mip:bin_limit} and \eqref{mip:bin_domain}. 
Last,
constraints \eqref{mip:number} and \eqref{mip:rank} are satisfied since $\bar{w}$ is regular. 

For the rest of the proof, 
we assume that the CGLP solution $\bar{w}$ is irregular. 

Suppose that the cut is \change{not strictly irregular}, 
and thus there is \change{an extended} regular CGLP solution $\tilde{w} = (\tilde{\alpha}, \tilde{\beta}, \{\tilde{u}^t, \tilde{v}^t\}_{t \in T})$ 
for which $(\bar{\alpha}, \bar{\beta}) = \lambda (\tilde{\alpha}, \tilde{\beta})$ for some $\lambda > 0$.
If $\lambda \leq 1$, 
then $(\lambda, \delta(\tilde{u}), \{ \tilde{u}^t, \tilde{v}^t \}_{t \in T})$ is feasible for RCV-MILP 
and thus $\check{\theta} \geq \lambda > 0$. 
If $\lambda > 1$, 
we can divide the same values by $\lambda$ to obtain another feasible CGLP solution 
$\hat{w} = \left(\frac{1}{\lambda} \tilde{\alpha}, \frac{1}{\lambda} \tilde{\beta},  \left\{ \frac{1}{\lambda} \tilde{u}^t, \frac{1}{\lambda}\tilde{v}^t \right\}_{t \in T} \right) = \left(\bar{\alpha}, \bar{\beta}, \left\{ \frac{1}{\lambda} \tilde{u}^t, \frac{1}{\lambda}\tilde{v}^t \right\}_{t \in T} \right)$. 
The construction of a RCV-MILP solution when $\lambda \leq 1$ applies to $\hat{w}$, 
thereby implying that 
$\left(1, \delta(\tilde{u}), \left\{ \frac{1}{\lambda} \tilde{u}^t, \frac{1}{\lambda} \tilde{v}^t \right\}_{t \in T} \right)$ is optimal for RCV-MILP and thus $\check{\theta} = 1$. 

Finally, 
suppose that the cut is \change{strictly} irregular, and thus there is no \change{extended} regular CGLP solution corresponding to cut $\bar{\alpha} x \geq \bar{\beta}$. 
Let us suppose, for contradiction, 
that there is a RCV-MILP solution $(\breve{\theta}, \breve{\delta}, \{ \breve{u}^t, \breve{v}^t \}_{t \in T})$ 
with $\breve{\theta} > 0$. 
Let \change{$\sigma := \sum_{t \in T} \sum_{i=1}^r \breve{u}^t_i + \sum_{t \in T} \sum_{i=1}^r \breve{v}^t_i$}. 
Since the cut separates $\bar{x}$ and $\breve{\theta} > 0$, 
then $(\bar{\alpha}, \bar{\beta}) \neq 0$ and thus $\sigma > 0$. 
Hence, 
if we divide the multipliers $\{ \breve{u}^t, \breve{v}^t \}_{t \in T}$ by $\sigma$,  
then normalization~\eqref{cglp:norm} is satisfied and $\left(\frac{\breve{\theta} }{\sigma} \bar{\alpha}, \frac{\breve{\theta} }{\sigma} \bar{\beta}, \left\{ \frac{1}{\sigma} \breve{u}^t, \frac{1}{\sigma} \breve{v}^t \right\}_{t \in T} \right)$ is a feasible \change{extended} regular CGLP solution: a contradiction. 
\endproof

\section{Numerical Procedure}\label{sec:algo}

Next we describe a numerical procedure to identify \change{strictly} irregular lift-and-project cuts, 
which addresses two issues with using RCV-MILP directly. 

First, 
finite numerical precision may lead to rounding errors and cause false negatives, 
thus overcounting the number of \change{strictly} irregular cuts. 
We address that by adding a relative tolerance parameter $\varepsilon$ on the coefficients of the cut. 
For example, 
if we want to determine if a cut  
$2 x_1 + 0.3 x_2 \geq 10$ is regular for $\varepsilon = 0.0001$, 
then we look for valid inequalities on each term where $\alpha_1 \in \theta [1.9998,2.0002]$, 
$\alpha_2 \in \theta [0.29997,0.30003]$, 
and $\beta \in \theta [9.999,10.001]$. 
\change{ We therefore avoid false negatives (often described as type II errors) at the price of tolerating false positives (often described as type I errors) when testing which cuts are regular. }  
This choice is intentional, 
since we are mainly interested in knowing which lift-and-project cuts are not regular.  
Furthermore, 
if $\varepsilon$ is sufficiently small, 
misclassifications are very unlikely.

Second, the number of subsets of $Q$ can be very large 
and many of those subsets might be irrelevant. 
We can address that by defining a set $\mathcal{Q}$ of subsets of $Q$ 
and then adding elements to this set as needed.

Hence, 
we define the 
Iterative RCV-MILP~(IRCV-MILP): 
\begin{align*}
\max ~\theta \\
& -\theta \varepsilon \leq \theta \bar{\alpha} - u^t \tilde{A} - v^t D^t \leq \theta \varepsilon, & t \in T  \\
& - \theta \varepsilon \leq \theta \bar{\beta} - u^t \tilde{b} - v^t d^t_0 \leq \theta \varepsilon, & t \in T \\
& \delta_j  - u^t_j  \geq  0, & j \in Q, t \in T \\
& \sum_{j \in Q} \delta_j  \leq  n \\
& \sum_{j \in N} \delta_j \leq \text{rank}(\tilde{A}_N), & N \subseteq \mathcal{Q}  \\
& u^t, v^t  \geq 0, & t \in T \\
& \delta_j \in \{0, 1\}, & j \in Q  \\ 
& \theta \in [0, 1] 
\end{align*}

Algorithm~\ref{alg:detect} uses IRCV-MILP to determine if a cut is regular, 
subject to false positives only. 
It finishes in finite time 
since each repetition of the loop corresponds to adding a different subset $N$ to $\mathcal{Q}$. 
In the unlikely worst case, 
Algorithm~\ref{alg:detect} terminates when all subsets of $Q$ have been added to $\mathcal{Q}$. 

\begin{algorithm}[h!]
\caption{Checks if there is a regular CGLP solution for a given cut}
\label{alg:detect}
{\AlgorithmFontSize
\begin{algorithmic}[1]
\Function{IsCutRegular}{$\bar{\alpha}, \bar{\beta}, \{ \bar{u}^t, \bar{v}^t \}_{t \in T}, ~ \varepsilon$}
\State $N \gets N(\bar{u})$
\If{rank$(\tilde{A}_N) = |N|$}
\State \Return True
\Comment{Original CGLP solution is regular}
\Else
\State $\mathcal{Q} \gets \emptyset$
\Loop
\State Get optimal solution $(\check{\theta}, \check{\delta}, \{ \check{u}^t, \check{v}^t \}_{t \in T})$ of IRCV-MILP
\State $N \gets N(\check{u})$
\If{$\bar{\theta}=0$}
\State \Return False
\Comment{There is no regular CGLP solution}
\ElsIf{rank$(\tilde{A}_N) < |N|$}
\State $\mathcal{Q} \gets \mathcal{Q} \cup \{ N \}$
\Comment{Loop has to be repeated}
\Else
\State \Return True
\Comment{Found regular CGLP solution}
\EndIf
\EndLoop
\EndIf
\EndFunction
\end{algorithmic}
}
\end{algorithm}

We can reduce the number of loop repetitions by preventing combinations of inequalities corresponding to parallel hyperplanes 
across different terms of the disjunction. 
For example, if rows $j_1$ and $j_2$ correspond to $x_i \geq 0$ and $x_i \leq 1$, 
then we can add the following inequality to IRCV-CGLP:

\begin{equation}\label{eq:bin_cut}
\delta_{j_1} + \delta_{j_2} \leq 1
\end{equation}

\section{Computational Experiments}\label{sec:exps}

We have run experiments to find \change{strictly} irregular lift-and-project cuts among the first round of cuts generated for 74 instances 
from the MIPLIB 2, 3, and 2003 benchmarks~\citep{MIPLIB, MIPLIB3, MIPLIB2003}. 
For each of those instances, 
we found an optimal solution $\bar{x}$ of the LP relaxation 
and generated a cut using the CGLP from each disjunction of the form~\eqref{eq:disj_k},
\change{i.e. $\bigvee_{K' \subseteq K} \{ x_k \geq 1 \text{ if } k \in K'; ~ x_k \leq 0 \text{ otherwise}\}$, where $0 < \bar{x}_k < 1$ for each $k \in K$,} 
with $|K|=2$ for all instances as well as $|K|=3$ 
and $|K|=4$ 
for smaller ones. 
Since the verification can be computationally expensive, 
the experiments are restricted to instances with at most 10,000 nonzeroes  
and each verification was interrupted if inconclusive after a predefined number of steps or time, 
which are both detailed later.
All code is written in C++ (gcc 4.8.2) and ran in Ubuntu 14.04.1 LTS on a machine with 
48 Intel(R) Core(TM) i7-4770 CPU @ 3.40GHz 
processors 
and
32 GB of RAM. 
No more than two copies were run in parallel and each was restricted to 10 GB of RAM. 
The formulations were solved with CPLEX version 12.6.3  
and matrix ranks were computed with Eigen\footnote{Available at \url{http://eigen.tuxfamily.org}}. 
Finally, 
we have run Algorithm~\ref{alg:detect} to determine if the cut is \change{strictly} irregular with $\varepsilon = 0.0001$ and 
constraints of the form~\eqref{eq:bin_cut} for the bounds of each variable $x_i, ~ i=1,\ldots, p$. 
\change{For brevity, we denote strictly irregular cuts as irregular cuts and the remaining cuts implied by the closure of regular cuts as regular cuts in the tables and figures.}

The most extensive experiments were run on 45 instances, 
where we counted the number of times that the loop of Algorithm~\ref{alg:detect} is repeated 
and set a time limit of 1 hour to interrupt the verification for each cut. 
We have chosen the subset of instances with at most 2,000 nonzeroes, 150 integer variables, or 50 rows; 
and we have also included larger instances from families of instances that nearly \change{ met those conditions}, 
namely \texttt{p0291}, \texttt{misc03}, \texttt{misc07}, and \texttt{pp08a}. 
For instances with at most 20 fractional values in $\bar{x}$ for integer variables, 
we generated and tested cuts with $|K|=2$ to $|K|=4$.
For instances with at most 30 fractional values, 
we generated and tested cuts with $|K|=2$ and $|K|=3$.
For other instances with at most 50 fractional values and \texttt{pp08a}, 
we only generated and tested cuts with $|K|=2$. 
\change{Some instances are discarded or only the results for less disjunctions are reported 
in cases where too many verifications timed out or the computer program ran out of resources.} 

\change{
We also verified if the cuts that are obtained with basic regular CGLP solutions map to split cuts. 
For a disjunction defined by $K$, 
we check if only one element in $v^t$ is positive for each $t \in T$ 
and if those positive multipliers correspond to inequalities on a same variable $x_k$, $k \in K$, 
in which case the convex set associated with the intersection cut is $\{ x: 0 \leq x_k \leq 1 \}$.}

The results for 36 instances having cuts with $|K|=2$ and $|K|=3$ are found in Table~\ref{tab:30}.
Additional results for the 22 instances having cuts with $|K|=4$ are found in Table~\ref{tab:20}.
We also compare the incidence of irregular CGLP solutions and cuts 
with $|K|=2$ to $|K|=4$ for those 22 instances in Table~\ref{tab:k234}.
The results for the remaining 9 instances having cuts with $|K|=2$ are found in Table~\ref{tab:50}. 
Since the number of cuts for each instance varies, 
we summarize a per-instance average of the percentages for each metric, 
which \change{ weighs} the results of each instance for equal contribution.

\begin{figure}[h!b]
\caption{Incidence of regular and irregular cuts from disjunctions with $|$K$|$ = 3 and $|$K$|$ = 4 according to the incidence of irregular cuts across all disjunctions from 2-subsets of K, i.e., all K' $\subset$ K such that $|$K'$|$=2.}
\label{fig:relation}
\centering
\vspace{2ex}
\begin{tikzpicture}[scale=0.55, font=\Large]
\begin{semilogyaxis}[
symbolic x coords={0, 1, 2, 3},
x tick label style={
/pgf/number format/1000 sep=},
xlabel={Irregular cuts from 2-subsets of $K$},
ylabel={Cuts with $|K|=3$},
enlarge x limits=0.2,
ymin=1,
ymax=80000,
legend style={/tikz/every even column/.append style={column sep=0.3cm}, at={(0.5,0.95)},
anchor=north,legend columns=-1},
ybar=0pt,
bar width=18pt,
height=200pt,
width=280pt,
xtick={0, 1, 2, 3},
]
    \addplot[ybar,fill=gray] coordinates {
        (0,   9464)
        (1,   1097)
        (2,   94)
        (3,   37)
    };
     \addplot[ybar,fill=black] coordinates {
        (0,   220)
        (1,   880)
        (2,   1918)
        (3,   5394)
    };
\legend{Regular, Irregular}
\end{semilogyaxis}
\end{tikzpicture}
~
\begin{tikzpicture}[scale=0.55, font=\Large]
\begin{semilogyaxis}[
symbolic x coords={0, 1, 2, 3, 4, 5, 6},
x tick label style={
/pgf/number format/1000 sep=},
xlabel={Irregular cuts from 2-subsets of $K$},
ylabel={Cuts with $|K|=4$},
enlarge x limits=0.1,
ymin=1,
ymax=80000,
legend style={/tikz/every even column/.append style={column sep=0.3cm}, at={(0.5,0.95)},
anchor=north,legend columns=-1},
ybar=0pt,
bar width=18pt,
height=200pt,
width=420pt,
xtick={0, 1, 2, 3, 4, 5, 6}
]
    \addplot[ybar,fill=gray] coordinates {
        (0,   7023)
        (1,   376)
        (2,   85)
        (3,   83)
        (4,   31)
        (5,   26)
        (6,   2)
    };
     \addplot[ybar,fill=black] coordinates {
        (0,   44)
        (1,   202)
        (2,   215)
        (3,   294)
        (4,   449)
        (5,   834)
        (6,   1008)
    };
\legend{Regular, Irregular}
\end{semilogyaxis}
\end{tikzpicture}
\end{figure}
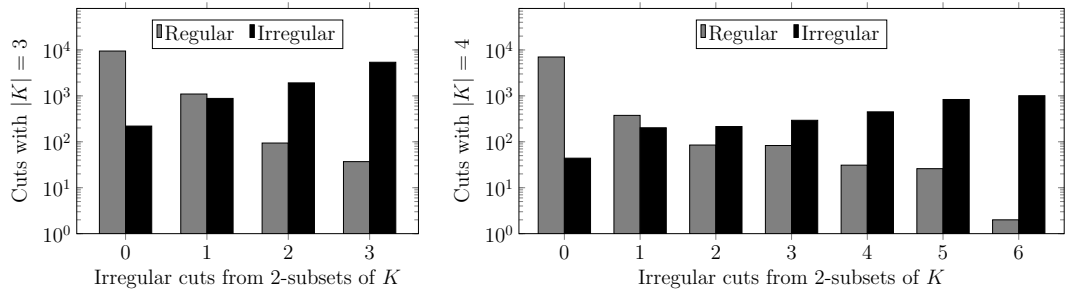

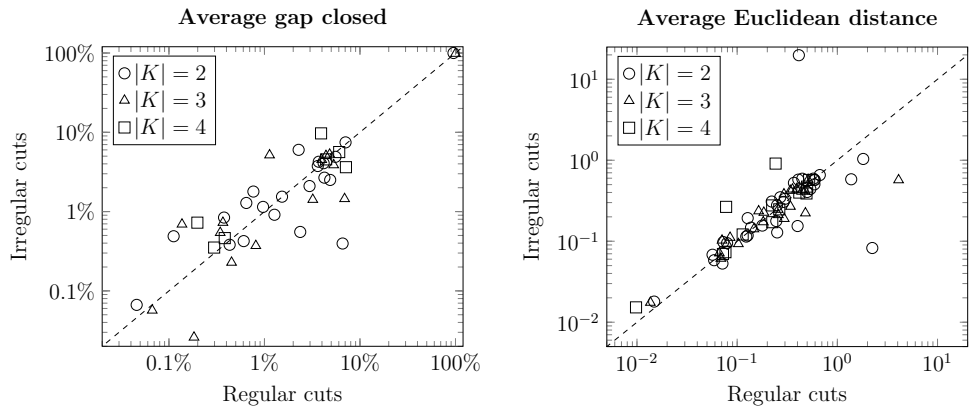
\begin{figure}
\caption{Comparison of regular vs. irregular cuts on instances and sizes of K where both types are found: 
the average gap closed is larger in 22 cases for irregular cuts vs. 24 cases for regular cuts; and  
the average Euclidean distance of the cut is larger in 37 cases for irregular cuts vs. 34 for regular cuts.}
\label{fig:average}
\vspace{2ex}
\centering
\begin{tikzpicture}[scale=0.7, font=\large]
\begin{loglogaxis}[
xmin=0.0002,xmax=1.2,ymin=0.0002,ymax=1.2,
legend pos=north west,
title={\textbf{Average gap closed}},
xlabel={Regular cuts},
ylabel={Irregular cuts},
xticklabels={},
yticklabels={},
extra x ticks={0.001,0.01,0.1,1.0},
extra x tick labels={0.1\%,1\%,10\%,100\%},
extra y ticks={0.001,0.01,0.1,1.0},
extra y tick labels={0.1\%,1\%,10\%,100\%},
mark options={scale=1.5},
scatter/classes={
    a={mark=o,draw=black},
    b={mark=triangle,draw=black},
    c={mark=square,draw=black}}]
\addplot[scatter,only marks,
    scatter src=explicit symbolic]
table[meta=label] {
x y label
0.0370359 0.0425187 a
0.00178059 0 a
0.0416333 0.0407672 a
0.00111949 0.00488765 a

0.00606061 0.00424244 a
0.00431954 0.00381552 a
0.00377524 0.00840736 a
0.0236855 0.00553718 a
0.000459208 0.000665816 a
0.0296344 0.0209893 a
0.0424191 0.0267853 a
0.95283 1 a
0.055375 0.0490549 a
0.0657962 0.0039568 a
0.00645651 0.0128453 a
0.0227886 0.0600803 a
0.00966934 0.0114583 a
0.0152724 0.0154133 a
0.0127213 0.00912104 a
0.00765658 0.0178954 a
0.0706127 0.0746072 a
0.0361581 0.0375662 a
0.0487495 0.0250707 a
0.0397199 0.0452326 b
0.00234955 0 b
0.0481266 0.0526873 b
0.00136921 0.00691565 b

0.00367356 0.0072725 b
0.00452396 0.00228297 b
0.00344215 0.00542116 b
0.00812051 0.00371094 b
0.000669163 0.000568381 b
0.0320474 0.0141711 b
0.0530922 0.0386986 b
1 0.984374 b
0.0440565 0.0505923 b
0.0689996 0.0144944 b
0.00182987 0.000259488 b
0.0113094 0.0515794 b
0.0442466 0.0455414 c

0.0605985 0.0563405 c

0.0391366 0.0969318 c
0.00198438 0.00725152 c

0.00297363 0.00351776 c
0.0038402 0.00462838 c

0.0712861 0.0362961 c

    };
\addplot[dashed, color=black] coordinates {(0.00001,0.00001)(2.0,2.0)};
\legend{$|K|=2$,$|K|=3$,$|K|=4$}
\end{loglogaxis}
\end{tikzpicture}
 ~ 
\begin{tikzpicture}[scale=0.7, font=\large]
\begin{loglogaxis}[
xmin=0.005,xmax=20,ymin=0.005,ymax=25,
legend pos=north west,
title={\textbf{Average Euclidean distance}},
xlabel={Regular cuts},
ylabel={Irregular cuts},
mark options={scale=1.5},
scatter/classes={
    a={mark=o,draw=black},
    b={mark=triangle,draw=black},
    c={mark=square,draw=black}}]
\addplot[scatter,only marks,
    scatter src=explicit symbolic]
table[meta=label] {
x y label
0.4081703125 0.575986724137931 a
0.0726919705882353 0.097724775 a
0.0589486818181818 0.0582902647058823 a
2.2331695240625 0.0824369846153846 a

0.178337109803921 0.156046333333333 a
0.127278528461538 0.117711052631579 a
0.245780538461538 0.17746 a
0.290118384615385 0.305883528301887 a
0.269354133333333 0.349043225609756 a
0.0567903375 0.0681068714285714 a
0.127378092307692 0.192848192 a
0.53788780952381 0.441488090909091 a

0.250923096153846 0.128432 a
0.22148750952381 0.308986545454545 a
0.014843685 0.01797165 a
0.5001935 0.40589625 a
0.412938222222222 19.81340258 a

0.071307 0.0529941333333333 a
0.590304181818182 0.5039906875 a

0.07933499 0.0948064636363636 a
0.137487833333333 0.147118407604167 a

0.401166444117647 0.153991885714286 a

0.12323532183908 0.114146 a
0.25019219882353 0.277843279166667 a
0.368799807142857 0.525536655555556 a
0.668139749425287 0.656129990624999 a
0.584180857142857 0.570229476545048 a
0.300771764150944 0.334486966779661 a
0.446748918167701 0.591507611111111 a
1.37539228471616 0.581930083333333 a
0.575720254084507 0.592452939393939 a
0.595665265478395 0.570916211739131 a
1.81258916045627 1.03694060808241 a
0.5109284 0.576051047058824 b
0.0792781606666665 0.0939184 b
0.0686498678571428 0.0612187107142857 b
0.0849902798666665 0.110990128270148 b

0.182196377777778 0.1727952 b
0.147089380848588 0.141515945119306 b
0.256812083333333 0.246121 b

0.339387571428571 0.421932751519757 b
0.0664882625 0.0712873857142857 b
0.184557529473684 0.226340730263158 b
0.477158838709677 0.415911417142857 b

0.295066210285714 0.189867307692308 b
0.257734375 0.211635428571429 b
0.340442114814815 0.268471753125 b
0.0136123333333333 0.0175109142857143 b
0.448423833333333 0.429915 b
4.08663406666667 0.57014775 b

0.102869 0.09280798 b
0.392670666666667 0.437218589341693 b

0.07064059 0.102617868 b
0.16388775 0.235750479822616 b

0.480090255263158 0.222476973076923 b

0.293173941573693 0.381689264543525 b
0.356226909043062 0.437796444444444 b
0.522902323943662 0.570935064748202 c

0.072153135483871 0.0695256153846154 c

0.447524588235294 0.445732 c
0.224104748329621 0.171823882352941 c

0.263257 0.247647 c

0.0764134 0.07319838 c
0.223652195583596 0.276818706864564 c
0.417035476190476 0.399070072093023 c
0.0775596926470588 0.265616 c

0.00978294 0.015251825 c
0.499483 0.472442 c
0.242012 0.909158444444445 c

0.113495 0.121273564705882 c

0.492007278899083 0.388749841176471 c

    };
\addplot[dashed, color=black] coordinates {(0.001,0.001)(50.0,50.0)};
\legend{$|K|=2$,$|K|=3$,$|K|=4$}
\end{loglogaxis}
\end{tikzpicture}
\end{figure}

Figure~\ref{fig:detection} shows the number of cuts identified as regular and irregular by Algorithm~\ref{alg:detect} according to the order of magnitude of repetitions of the loop upon termination. 
Figure~\ref{fig:relation} compares the regularity of cuts from each $K$ of size 3 and 4
with that of cuts from subsets of $K$ of size 2 
among the 97\% of the cases where no cut for a subset timed out. 
\change{In other words, for each set $K = \{k_1, \ldots, k_m\}, m \in \{3, 4\},$ that defines a simple $m$-branch split disjunction from which we obtain a regular or irregular cut, 
we count how many irregular cuts are obtained from each subset $K' \subset K, |K'| = 2,$ which defines a simple $2$-branch split disjunction. There are 3 of those subsets of size 2 in $K$ if $|K|=3$ and 6 if $|K|=4$.} 
Figure~\ref{fig:average} compares the average gap closed and the average Euclidean distance of the cuts generated with respect to the fractional solution $\bar{x}$ on all disjunction sizes where both types of cuts are observed on each of the 45 instances. 
\change{Figure~\ref{fig:total_gap} shows the total gap closed with and without irregular cuts in the cases in which both types of cuts are observed.}

\begin{figure}[h!]
\caption{Incidence of cuts assessed by Algorithm~\ref{alg:detect} by the order of repetitions of the loop upon termination.}
\label{fig:detection}
\centering
\vspace{2ex}
\begin{tikzpicture}[scale=0.5, font=\Large]
\begin{axis}[
symbolic x coords={0, 1, 4, 16, 64, 256, 1024, 4096},
x tick label style={
/pgf/number format/1000 sep=},
xlabel={Loop repetitions},
ylabel={Cuts with $|K|=2$},
enlarge x limits=0.1,
ybar stacked, 
ymin=0,
ymax=11950,
scaled y ticks=false,
bar width=14pt,
legend style={/tikz/every even column/.append style={column sep=0.3cm}, at={(0.5,0.95)},
anchor=north, legend columns=-1},
xtick={0,1,4,64,1024}
]
    \addplot[fill=gray] coordinates {
(0, 5055)
(1, 438)
(4, 307)
(16, 230)
(64, 234)
(256, 197)
(1024, 101)
(4096, 22)
    };
     \addplot[fill=black] coordinates {
(0, 0)
(1, 3635)
(4, 237)
(16, 34)
(64, 33)
(256, 24)
(1024, 17)
(4096, 6)
    };
\legend{Regular,Irregular}
\end{axis}
\end{tikzpicture}
\begin{tikzpicture}[scale=0.5, font=\Large]
\begin{axis}[
symbolic x coords={0, 1, 4, 16, 64, 256, 1024, 4096},
x tick label style={
/pgf/number format/1000 sep=},
xlabel={Loop repetitions},
ylabel={Cuts with $|K|=3$},
enlarge x limits=0.1,
ybar stacked, 
ymin=0,
ymax=11950,
scaled y ticks=false,
bar width=14pt,
legend style={/tikz/every even column/.append style={column sep=0.3cm}, at={(0.5,0.95)},
anchor=north, legend columns=-1},
xtick={0,1,4,64,1024}
]
    \addplot[ybar,fill=gray] coordinates {
(0, 9497)
(1, 895)
(4, 280)
(16, 144)
(64, 172)
(256, 123)
(1024, 63)
(4096, 2)
    };
     \addplot[ybar,fill=black] coordinates {
(0, 0)
(1, 7471)
(4, 851)
(16, 81)
(64, 107)
(256, 70)
(1024, 13)
(4096, 5)
    };
\legend{Regular,Irregular}
\end{axis}
\end{tikzpicture}
\begin{tikzpicture}[scale=0.5, font=\Large]
\begin{axis}[
symbolic x coords={0, 1, 4, 16, 64, 256, 1024, 4096},
x tick label style={
/pgf/number format/1000 sep=},
xlabel={Loop repetitions},
ylabel={Cuts with $|K|=4$},
enlarge x limits=0.1,
ybar stacked, 
ymin=0,
ymax=11950,
scaled y ticks=false,
bar width=14pt,
legend style={/tikz/every even column/.append style={column sep=0.3cm}, at={(0.5,0.95)},
anchor=north, legend columns=-1},
xtick={0,1,4,64,1024}
]
    \addplot[ybar,fill=gray] coordinates {
(0, 6834)
(1, 401)
(4, 272)
(16, 112)
(64, 95)
(256, 93)
(1024, 42)
(4096, 4)
    };
     \addplot[ybar,fill=black] coordinates {
(0, 0)
(1, 2891)
(4, 142)
(16, 26)
(64, 27)
(256, 22)
(1024, 22)
(4096, 15)
    };
\legend{Regular,Irregular}
\end{axis}
\end{tikzpicture}
\end{figure}
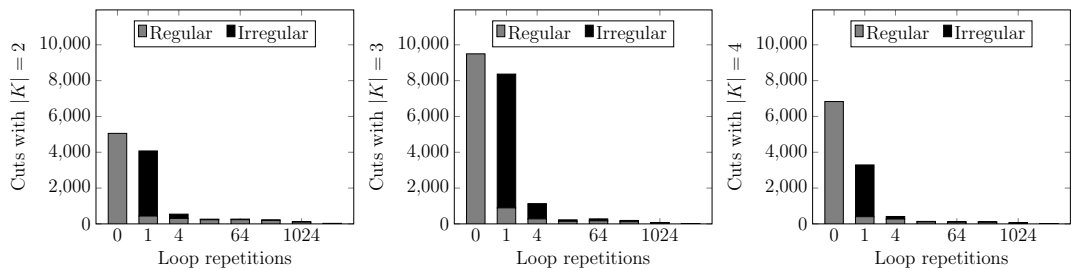

\newgeometry{left=0.5in,right=0.5in,top=0.5in,bottom=0.5in}
\begin{landscape}
\begin{table}
\footnotesize
\caption{Results for instances where lift-and-project cuts with $|$K$|$=2 and $|$K$|$=3 are generated \change{ and} tested. For each instance and size of K, 
we report the regularity of the CGLP solution, \change{when those solutions define split cuts,} and the regularity of the cut by running Algorithm~\ref{alg:detect} with time limit of one hour.}
\label{tab:30}
\centering
\vspace{2ex}
\begin{tabular}{@{\extracolsep{4pt}}cccccccccccccc}
&& \multicolumn{6}{c}{Lift-and-project cuts with $|K|=2$} & \multicolumn{6}{c}{Lift-and-project cuts with $|K|=3$} \\
\cline{3-8}
\cline{9-14}
\noalign{\vskip2.5pt}
& Fractional & \multicolumn{2}{c}{CGLP basis} & \multicolumn{4}{c}{Cut} & \multicolumn{2}{c}{CGLP basis} & \multicolumn{4}{c}{Cut} \\
\cline{3-4}
\cline{5-8}
\cline{9-10}
\cline{11-14}
\noalign{\vskip2.5pt}
Instance & Variables & Regular & Irregular & \change{Split} & Regular & Irregular & Unknown & Regular & Irregular & \change{Split} & Regular & Irregular & Unknown \\
\cline{1-1}
\cline{2-2}
\cline{3-4}
\cline{5-5}
\cline{6-8}
\cline{9-10}
\cline{11-11}
\cline{12-14}
\noalign{\vskip2.5pt}
air01 & 5 & 4 & 6 & \change{0} & 4 & 6 & 0 & 1 & 9 & \change{0} & 1 & 5 & 4 \\ 
bell5 & 25 & 240 & 60 & \change{64} & 280 & 9 & 11 & 1787 & 513 & \change{82} & 2090 & 63 & 147 \\ 
blend2 & 6 & 4 & 11 & \change{2} & 9 & 5 & 1 & 8 & 12 & \change{1} & 15 & 4 & 1 \\ 
bm23 & 6 & 0 & 15 & \change{0} & 0 & 15 & 0 & 0 & 20 & \change{0} & 0 & 20 & 0 \\ 
enigma & 8 & 4 & 24 & \change{2} & 25 & 0 & 3 & 4 & 52 & \change{3} & 45 & 0 & 11 \\ 
flugpl & 10 & 16 & 29 & \change{5} & 16 & 29 & 0 & 33 & 87 & \change{1} & 35 & 85 & 0 \\ 
gt2 & 11 & 55 & 0 & \change{0} & 55 & 0 & 0 & 153 & 12 & \change{0} & 164 & 0 & 1 \\ 
khb05250 & 19 & 10 & 161 & \change{10} & 11 & 160 & 0 & 0 & 969 & \change{0} & 6 & 957 & 6 \\ 
lseu & 11 & 43 & 12 & \change{31} & 51 & 4 & 0 & 92 & 73 & \change{46} & 150 & 15 & 0 \\ 
markshare1 & 6 & 0 & 15 & \change{0} & 0 & 15 & 0 & 0 & 20 & \change{0} & 0 & 20 & 0 \\ 
markshare2 & 7 & 10 & 11 & \change{10} & 10 & 11 & 0 & 9 & 26 & \change{9} & 10 & 25 & 0 \\ 
mas74 & 12 & 65 & 1 & \change{61} & 66 & 0 & 0 & 203 & 17 & \change{195} & 220 & 0 & 0 \\ 
mas76 & 11 & 55 & 0 & \change{55} & 55 & 0 & 0 & 165 & 0 & \change{165} & 165 & 0 & 0 \\ 
misc01 & 12 & 64 & 2 & \change{47} & 64 & 0 & 2 & 166 & 54 & \change{66} & 175 & 13 & 32 \\ 
misc02 & 8 & 26 & 2 & \change{19} & 26 & 2 & 0 & 32 & 24 & \change{15} & 32 & 14 & 10 \\ 
misc03 & 14 & 87 & 4 & \change{76} & 87 & 1 & 3 & 334 & 30 & \change{256} & 334 & 0 & 30 \\ 
misc05 & 11 & 17 & 38 & \change{14} & 21 & 22 & 12 & 18 & 147 & \change{6} & 27 & 64 & 74 \\ 
misc07 & 16 & 115 & 5 & \change{108} & 115 & 0 & 5 & 502 & 58 & \change{442} & 502 & 0 & 58 \\ 
mod008 & 5 & 0 & 10 & \change{0} & 4 & 6 & 0 & 0 & 10 & \change{0} & 3 & 7 & 0 \\ 
mod013 & 5 & 4 & 6 & \change{4} & 6 & 4 & 0 & 2 & 8 & \change{1} & 6 & 3 & 1 \\ 
modglob & 30 & 2 & 433 & \change{1} & 32 & 403 & 0 & 0 & 4060 & \change{0} & 75 & 3983 & 2 \\ 
p0033 & 6 & 13 & 2 & \change{11} & 15 & 0 & 0 & 19 & 1 & \change{13} & 20 & 0 & 0 \\ 
p0040 & 4 & 5 & 1 & \change{2} & 6 & 0 & 0 & 3 & 1 & \change{0} & 4 & 0 & 0 \\ 
p0201 & 20 & 189 & 1 & \change{144} & 190 & 0 & 0 & 1124 & 16 & \change{448} & 1140 & 0 & 0 \\ 
p0282 & 26 & 241 & 84 & \change{185} & 255 & 48 & 22 & 1585 & 1015 & \change{983} & 1703 & 471 & 426 \\ 
p0291 & 10 & 30 & 15 & \change{16} & 34 & 7 & 4 & 70 & 50 & \change{13} & 76 & 26 & 18 \\ 
pipex & 6 & 3 & 12 & \change{0} & 8 & 7 & 0 & 0 & 20 & \change{0} & 8 & 7 & 5 \\ 
pk1 & 15 & 1 & 104 & \change{1} & 9 & 96 & 0 & 0 & 455 & \change{0} & 4 & 451 & 0 \\ 
rgn & 14 & 40 & 51 & \change{39} & 65 & 25 & 1 & 76 & 288 & \change{70} & 190 & 152 & 22 \\ 
sample2 & 12 & 18 & 48 & \change{7} & 21 & 44 & 1 & 28 & 192 & \change{0} & 31 & 175 & 14 \\ 
sentoy & 8 & 1 & 27 & \change{1} & 11 & 17 & 0 & 0 & 56 & \change{0} & 28 & 28 & 0 \\ 
stein9 & 6 & 9 & 6 & \change{1} & 13 & 2 & 0 & 2 & 18 & \change{0} & 12 & 8 & 0 \\ 
stein15 & 12 & 9 & 57 & \change{1} & 13 & 53 & 0 & 0 & 220 & \change{0} & 0 & 220 & 0 \\ 
stein27 & 21 & 42 & 168 & \change{15} & 45 & 164 & 1 & 8 & 1322 & \change{1} & 14 & 1316 & 0 \\ 
vpm1 & 15 & 100 & 5 & \change{98} & 102 & 3 & 0 & 441 & 14 & \change{436} & 450 & 5 & 0 \\ 
vpm2 & 30 & 314 & 121 & \change{311} & 390 & 38 & 7 & 2632 & 1428 & \change{2615} & 3441 & 461 & 158 \\ 
\cline{1-2}
\cline{3-4}
\cline{5-5}
\cline{6-8}
\cline{9-10}
\cline{11-11}
\cline{12-14}
\firstcline{1-2}
\noalign{\vskip1.5pt} 
\multicolumn{2}{c}{Instance average} & 51.3\% & 48.7\% & \change{34.6\%} & 62.9\% & 35.0\% & 2.0\% & 40.8\% & 59.2\% & \change{22.5\%} &  54.5\% & 38.5\% & 7.0\% \\
\end{tabular}
\end{table}
\end{landscape}
\restoregeometry

\begin{figure}
\caption{\change{Comparison of total gap closed with and without irregular cuts on instances and sizes of $K$ where both types are found: the total gap closed with irregular cuts is larger in 27 cases.}}
\label{fig:total_gap}
\vspace{2ex}
\centering
\begin{tikzpicture}[scale=0.7, font=\large]
\begin{loglogaxis}[
xmin=0.0002,xmax=1.2,ymin=0.0002,ymax=1.2,
legend pos=south east,
title={\textbf{Total gap closed}},
xlabel={Without irregular cuts},
ylabel={With irregular cuts},
xticklabels={},
yticklabels={},
extra x ticks={0.001,0.01,0.1,1.0},
extra x tick labels={0.1\%,1\%,10\%,100\%},
extra y ticks={0.001,0.01,0.1,1.0},
extra y tick labels={0.1\%,1\%,10\%,100\%},
mark options={scale=1.5},
scatter/classes={
    a={mark=o,draw=black},
    b={mark=triangle,draw=black},
    c={mark=square,draw=black}}]
\addplot[scatter,only marks,
    scatter src=explicit symbolic]
table[meta=label] {
x y label
0.117413 0.122783 a
0.0420584 0.0420584 a
0.0815434 0.0870955 a
0.0207465 0.0525415 a

0.101819 0.101819 a
0.128828 0.141364 a
0.00787746 0.0215218 a
0.117188 0.117188 a
0.00073311 0.000845424 a
0.0441177 0.0441177 a
0.159358 0.159358 a
1 1 a
0.594643 0.729448 a
0.533948 0.533948 a
0.486185 0.508221 a
0.85374 0.85374 a
0.37631 0.461977 a
0.254511 0.528676 a
0.348186 0.365088 a
0.286943 0.286943 a
0.662571 0.662572 a
0.48288 0.647158 a
0.632441 0.63632 a

0.123615 0.125421 b
0.0420591 0.0420591 b
0.0805695 0.0946617 b
0.0428249 0.0653727 b

0.101819 0.101819 b
0.117896 0.179527 b
0.00769805 0.0200051 b
0.117188 0.117188 b
0.00110316 0.00110316 b
0.0441178 0.0441178 b
0.169612 0.363585 b
1 1 b
0.182035 0.692912 b
0.409876 0.409876 b
0.130656 0.15178 b
0.285575 0.285575 b

0.0984663 0.122643 c

0.103575 0.103575 c

0.448608 0.448608 c
0.0836376 0.0872745 c

0.00300512 0.00793078 c
0.080644 0.104167 c

0.159363 0.230303 c

    };
\addplot[dashed, color=black] coordinates {(0.00001,0.00001)(2.0,2.0)};
\legend{$|K|=2$,$|K|=3$,$|K|=4$}
\end{loglogaxis}
\end{tikzpicture}
\end{figure}

\begin{table}
\footnotesize
\caption{Additional results for instances where lift-and-project cuts with $|$K$|$=4 are also generated and tested.}
\label{tab:20}
\centering
\vspace{2ex}
\begin{tabular}{@{\extracolsep{4pt}}cccccccc}
&& \multicolumn{6}{c}{Lift-and-project cuts with $|K|=4$} \\
\cline{3-8}
\noalign{\vskip2.5pt}
& Fractional & \multicolumn{2}{c}{CGLP basis} & \multicolumn{4}{c}{Cut}  \\
\cline{3-4}
\cline{5-8}
\noalign{\vskip2.5pt}
Instance & Variables & Regular & Irregular & Split & Regular & Irregular & Unknown \\
\cline{1-1}
\cline{2-2}
\cline{3-4}
\cline{5-5}
\cline{6-8}
\noalign{\vskip2.5pt}
blend2 & 6 & 4 & 11 & \change{0} & 6 & 9 & 0 \\ 
bm23 & 6 & 0 & 15 & \change{0} & 0 & 15 & 0 \\ 
enigma & 8 & 4 & 66 & \change{3} & 68 & 1 & 1 \\ 
flugpl & 10 & 52 & 158 & \change{0} & 71 & 139 & 0 \\ 
gt2 & 11 & 179 & 151 & \change{0} & 238 & 1 & 91 \\ 
markshare1 & 6 & 0 & 15 & \change{0} & 0 & 15 & 0 \\ 
markshare2 & 7 & 1 & 34 & \change{1} & 1 & 34 & 0 \\ 
misc03 & 14 & 850 & 151 & \change{508} & 856 & 0 & 145 \\ 
mod008 & 5 & 0 & 5 & \change{0} & 1 & 4 & 0 \\ 
mod013 & 5 & 2 & 3 & \change{0} & 3 & 2 & 0 \\ 
p0033 & 6 & 15 & 0 & \change{8} & 15 & 0 & 0 \\ 
p0040 & 4 & 0 & 1 & \change{0} & 1 & 0 & 0 \\ 
p0201 & 20 & 4191 & 654 & \change{560} & 4753 & 0 & 92 \\ 
p0291 & 10 & 106 & 104 & \change{6} & 109 & 34 & 67 \\ 
pipex & 6 & 0 & 15 & \change{0} & 3 & 5 & 7 \\ 
pk1 & 15 & 0 & 1365 & \change{0} & 0 & 1365 & 0 \\ 
rgn & 14 & 74 & 927 & \change{62} & 317 & 539 & 145 \\ 
sample2 & 12 & 17 & 478 & \change{0} & 21 & 430 & 44 \\ 
sentoy & 8 & 0 & 70 & \change{0} & 31 & 39 & 0 \\ 
stein9 & 6 & 0 & 15 & \change{0} & 12 & 3 & 0 \\ 
stein15 & 12 & 0 & 495 & \change{0} & 0 & 495 & 0 \\ 
vpm1 & 15 & 1339 & 26 & \change{1336} & 1347 & 17 & 1 \\ 
\cline{1-2}
\cline{3-4}
\cline{5-5}
\cline{6-8}
\firstcline{1-2}
\noalign{\vskip1.5pt} 
\multicolumn{2}{c}{Instance average} & 26.6\% & 73.4\% & \change{5.2\%} & 47.3\% & 46\% & 6.7\% \\
\cline{1-8}
\firstcline{1-8}
\noalign{\vskip1.5pt} 
\multicolumn{2}{c}{\change{Average with $|K|=3$}} & \change{35.1\%} & \change{64.9\%} & \change{15.8\%} & \change{53.3\%} & \change{42.3\%} & \change{4.4\%} \\
\multicolumn{2}{c}{\change{Average with $|K|=2$}} & \change{43.7\%} & \change{56.3\%} & \change{26.4\%} & \change{59.6\%} & \change{38.9\%} & \change{1.5\%} \\
\end{tabular}
\end{table}

\begin{table}
\footnotesize
\caption{Frequency of irregular bases and cuts with $|$K$|$=2, 3, and 4.}
\label{tab:k234}
\centering
\vspace{2ex}
\begin{tabular}{@{\extracolsep{4pt}}ccccccc}
\noalign{\vskip2.5pt}
& \multicolumn{2}{c}{Irregularity with $|K| = 2$} & \multicolumn{2}{c}{Irregularity with $|K| = 3$} & \multicolumn{2}{c}{Irregularity with $|K| = 4$}  \\
\cline{2-3}
\cline{4-5}
\cline{6-7}
\noalign{\vskip2.5pt}
Instance & Basis & Cut & Basis & Cut & Basis & Cut \\
\cline{1-1}
\cline{2-2}
\cline{3-3}
\cline{4-4}
\cline{5-5}
\cline{6-6}
\cline{7-7}
\noalign{\vskip2.5pt}
blend2 & 73.3\% & 33.3\% & 60\% & 20\% & 73.3\% & 60\%  \\ 
bm23 & 100\% & 100\% & 100\% & 100\% & 100\% & 100\%  \\ 
enigma & 85.7\% & 0\% & 92.9\% & 0\% & 94.3\% & 1.4\%  \\ 
flugpl & 64.4\% & 64.4\% & 72.5\% & 70.8\% & 75.2\% & 66.2\%  \\ 
gt2 & 0\% & 0\% & 7.3\% & 0\% & 45.8\% & 0.3\%  \\ 
markshare1 & 100\% & 100\% & 100\% & 100\% & 100\% & 100\%  \\ 
markshare2 & 52.4\% & 52.4\% & 74.3\% & 71.4\% & 97.1\% & 97.1\%  \\ 
misc03 & 4.4\% & 1.1\% & 8.2\% & 0\% & 15.1\% & 0\%  \\ 
mod008 & 100\% & 60\% & 100\% & 70\% & 100\% & 80\%  \\ 
mod013 & 60\% & 40\% & 80\% & 30\% & 60\% & 40\%  \\ 
p0033 & 13.3\% & 0\% & 5\% & 0\% & 0\% & 0\%  \\ 
p0040 & 16.7\% & 0\% & 25\% & 0\% & 100\% & 0\%  \\ 
p0201 & 0.5\% & 0\% & 1.4\% & 0\% & 13.5\% & 0\%  \\ 
p0291 & 33.3\% & 15.6\% & 41.7\% & 21.7\% & 49.5\% & 16.2\%  \\ 
pipex & 80\% & 46.7\% & 100\% & 35\% & 100\% & 33.3\%  \\ 
pk1 & 99\% & 91.4\% & 100\% & 99.1\% & 100\% & 100\%  \\ 
rgn & 56\% & 27.5\% & 79.1\% & 41.8\% & 92.6\% & 53.8\%  \\ 
sample2 & 72.7\% & 66.7\% & 87.3\% & 79.5\% & 96.6\% & 86.9\%  \\ 
sentoy & 96.4\% & 60.7\% & 100\% & 50\% & 100\% & 55.7\%  \\ 
stein9 & 40\% & 13.3\% & 90\% & 40\% & 100\% & 20\%  \\ 
stein15 & 86.4\% & 80.3\% & 100\% & 100\% & 100\% & 100\%  \\ 
vpm1 & 4.8\% & 2.9\% & 3.1\% & 1.1\% & 1.9\% & 1.2\%  \\ 
\cline{1-1}
\cline{2-2}
\cline{3-3}
\cline{4-4}
\cline{5-5}
\cline{6-6}
\cline{7-7}
\firstcline{1-1}
\noalign{\vskip1.5pt} 
Instance average & 56.3\% & 38.9\% & 64.9\% & 42.3\% & 73.4\% & 46\%  \\ 
\end{tabular}
\end{table}

\begin{table}
\footnotesize
\caption{Results for instances where only lift-and-project cuts with $|K|=2$ are generated and tested.}
\label{tab:50}
\centering
\vspace{2ex}
\begin{tabular}{@{\extracolsep{4pt}}cccccccc}
&& \multicolumn{6}{c}{Lift-and-project cuts with $|K|=2$} \\
\cline{3-8}
\noalign{\vskip2.5pt}
& Fractional & \multicolumn{2}{c}{CGLP basis} & \multicolumn{4}{c}{Cut}  \\
\cline{3-4}
\cline{5-8}
\noalign{\vskip2.5pt}
Instance & Variables & Regular & Irregular & Split & Regular & Irregular & Unknown \\
\cline{1-1}
\cline{2-2}
\cline{3-4}
\cline{5-5}
\cline{6-8}
\noalign{\vskip2.5pt}
bell3a & 32 & 404 & 92 & \change{26} & 458 & 12 & 26 \\
bell3b & 36 & 442 & 188 & \change{124} & 568 & 33 & 29 \\
bell4 & 46 & 887 & 148 & \change{289} & 966 & 18 & 51 \\
dcmulti & 49 & 503 & 673 & \change{439} & 526 & 631 & 19 \\
egout & 40 & 16 & 764 & \change{0} & 435 & 320 & 25 \\
noswot & 22 & 189 & 42 & \change{90} & 199 & 0 & 32 \\
p0548 & 48 & 313 & 815 & \change{145} & 648 & 138 & 342 \\
pp08a & 53 & 35 & 1343 & \change{0} & 35 & 1343 & 0 \\
pp08aCUTS & 46 & 430 & 605 & \change{0} & 636 & 295 & 104 \\
\cline{1-2}
\cline{3-4}
\cline{5-5}
\cline{6-8}
\firstcline{1-2}
\noalign{\vskip1.5pt} 
\multicolumn{2}{c}{Instance average} & 48.4\% & 51.6\% & \change{15.8\%} & 64.9\% & 26.9\% & 8.2\% \\
\end{tabular}
\end{table}

For the remaining 29 instances reported in Table~\ref{tab:more}, 
we repeat the loop in Algorithm~\ref{alg:detect} once. 
This set includes some instances previously excluded.

\begin{table}[h!]
\footnotesize
\caption{Remaining results for larger instances where lift-and-project cuts with $|K|=2$ are generated and tested, but the loop of Algorithm~\ref{alg:detect} is only repeated once.}
\label{tab:more}
\centering
\vspace{2ex}
\begin{tabular}{@{\extracolsep{4pt}}cccccccc}
&& \multicolumn{6}{c}{Lift-and-project cuts with $|K|=2$} \\
\cline{3-8}
\noalign{\vskip2.5pt}
& Fractional & \multicolumn{2}{c}{CGLP basis} & \multicolumn{4}{c}{Cut}  \\
\cline{3-4}
\cline{5-8}
\noalign{\vskip2.5pt}
Instance & Variables & Regular & Irregular & Split & Regular & Irregular & Unknown \\
\cline{1-1}
\cline{2-2}
\cline{3-4}
\cline{5-5}
\cline{6-8}
\noalign{\vskip2.5pt}
aflow30a & 31 & 155 & 310 & \change{130} & 157 & 296 & 12 \\ 
aflow40b & 38 & 347 & 356 & \change{337} & 349 & 240 & 114 \\ 
cracpb1 & 40 & 389 & 391 & \change{364} & 408 & 259 & 113 \\ 
dsbmip & 48 & 464 & 664 & \change{453} & 871 & 36 & 221 \\ 
fiber & 42 & 687 & 174 & \change{657} & 713 & 42 & 106 \\ 
fixnet3 & 69 & 198 & 2148 & \change{0} & 228 & 0 & 2118 \\ 
fixnet4 & 67 & 180 & 2031 & \change{0} & 196 & 0 & 2015 \\ 
fixnet6 & 60 & 139 & 1631 & \change{0} & 153 & 0 & 1617 \\ 
gen & 41 & 505 & 315 & \change{483} & 511 & 58 & 251 \\ 
gesa2 & 58 & 704 & 949 & \change{655} & 749 & 7 & 897 \\ 
gesa2\_o & 73 & 1356 & 1272 & \change{1251} & 1383 & 763 & 482 \\ 
gesa3 & 85 & 1785 & 1785 & \change{1717} & 1820 & 589 & 1161 \\ 
gesa3\_o & 100 & 2175 & 2775 & \change{2092} & 2266 & 1555 & 1129 \\ 
glass4 & 72 & 2556 & 0 & \change{2556} & 2556 & 0 & 0 \\ 
harp2 & 30 & 435 & 0 & \change{0} & 435 & 0 & 0 \\ 
l152lav & 51 & 254 & 1021 & \change{228} & 292 & 904 & 79 \\ 
lp4l & 23 & 50 & 203 & \change{27} & 73 & 121 & 59 \\ 
opt1217 & 27 & 74 & 277 & \change{48} & 88 & 106 & 157 \\ 
p2756 & 77 & 1521 & 1405 & \change{452} & 2190 & 282 & 454 \\ 
qiu & 36 & 474 & 156 & \change{404} & 474 & 16 & 140 \\ 
qnet1 & 47 & 578 & 503 & \change{564} & 645 & 12 & 424 \\ 
qnet1\_o & 11 & 20 & 35 & \change{19} & 32 & 0 & 23 \\ 
rout & 35 & 304 & 291 & \change{298} & 308 & 220 & 67 \\ 
set1al & 218 & 239 & 23414 & \change{19} & 276 & 20940 & 2437 \\ 
set1ch & 138 & 122 & 9331 & \change{30} & 150 & 8786 & 517 \\ 
set1cl & 220 & 180 & 23910 & \change{19} & 225 & 21393 & 2472 \\ 
timtab1 & 134 & 2829 & 6082 & \change{2360} & 2942 & 5935 & 34 \\ 
timtab2 & 236 & 8929 & 18801 & \change{6856} & 9172 & 18027 & 531 \\ 
tr12-30 & 348 & 47 & 60331 & \change{0} & 605 & 59046 & 727 \\ 
\cline{1-2}
\cline{3-4}
\cline{5-5}
\cline{6-8}
\firstcline{1-2}
\noalign{\vskip1.5pt} 
\multicolumn{2}{c}{Instance average} &38.7\% & 61.3\% & \change{30.5\%} & 43.0\% & 31.8\% & 25.2\% \\
\end{tabular}
\end{table}

\section{Discussion}\label{sec:disc}

The experiments above evidence some trends on irregular CGLP bases and \change{strictly irregular} cuts, 
from which we can draw observations about their incidence \change{and relation to split cuts}.

First, 
there is a substantial difference between the number of irregular CGLP bases and that of \change{strictly} irregular cuts. 
In all tables, 
the per-instance averages differ in at least 10\% with respect to the total number of cuts.
Therefore, an irregular CGLP basis often does not guarantee that the corresponding cut is \change{strictly} irregular.

Second, most regular and irregular cuts can nevertheless be identified with the first test applied on them. 
Such paradox is in part explained by the fact that \change{strictly} irregular cuts are less frequent than regular cuts in our experiments. 
Figure~\ref{fig:detection} shows that most regular cuts are identified immediately because they come from regular CGLP solutions, 
whereas most irregular cuts are identified in the first repetition of the loop of Algorithm~\ref{alg:detect}. 
The same is true for the larger problems reported in Table~\ref{tab:more}: 
with a few exceptions,  
the majority of the cuts is identified as regular or irregular even though the loop is never repeated. 
In other words, 
restricting the number of nonzero multipliers of the linear relaxation to $n$ across all terms of the CGLP formulation 
is an effective way of determining if a cut is irregular. 
Conversely, 
most \change{strictly} irregular cuts are those that can only be derived by using more than $n$ rows of the linear relaxation. 
Hence,  
the linear dependence among 
at most $n$ 
rows is a secondary factor for cut irregularity in practice.

Third, 
simple $t$-branch split disjunctions with more terms yield more irregular cuts 
and even more irregular CGLP bases. 
We can observe that first from the difference in per-instance averages of the corresponding columns 
for $|K|=2$ and $|K|=3$ in Table~\ref{tab:30} 
and from $|K|=2$ to $|K|=4$ in Table~\ref{tab:k234}.
In addition,  
the frequency of irregular CGLP bases and cuts for each instance often increases with the size of $K$, 
the former more than the latter. 
For the 36 instance reported in Table~\ref{tab:30}, 
the frequency of irregular CGLP bases for $|K|=3$ is higher in 29 instances and lower in 3 when compared to the case of $|K|=2$, 
whereas the frequency of irregular cuts is higher in 17 and lower in 9. 
For the 22 instances reported in Table~\ref{tab:k234}, 
there are proportionally more irregular CGLP bases for $|K|=4$ in 15 instances and less in 2 instances when compared to $|K|=2$, 
while the ratio of irregular cuts is higher in 12 instances and lower in 4 instances.
\change{Similarly, disjunction with more terms yield proportonially less basic regular CGLP solutions mapping to split cuts.} 
Hence, a disjunction with more terms seems more likely to yield lift-and-project cuts 
that are not 
obtainable 
from a split disjunction.

Fourth, 
some disjunctions are more likely to yield irregular cuts than others. 
\change{
For simple $t$-branch split disjunctions, 
we observe in Figure~\ref{fig:relation} that an irregular cut obtained from a disjunction defined by a set $K$ is such that 
the disjunctions defined by subsets of $K$ of size 2 often yield irregular cuts as well. 
Conversely, if we know which sets of size 2 define disjunctions yielding irregular cuts, then a union of such sets defines a disjunction on more terms that is more likely to yield irregular cuts.} 
For example, 
a disjunction on variables $x_A$, $x_B$, and $x_C$ yields an irregular cut with increasing probability 
as more irregular cuts are obtained among disjunctions on $x_A$ and $x_B$, $x_A$ and $x_C$, and $x_B$ and $x_C$. 
\change{The probability of obtaining an irregular cut from a disjunction defined by a set $K$ of size 3 or 4 is already close to 50\% if at least one subset of $K$ of size 2 defines a disjunction yielding an irregular cut.}  
A simpler hypothesis that could be ventured is that 
the presence of certain variables in the disjunction 
makes it more likely that the CGLP would yield irregular cuts, 
although it is yet to be determined what makes disjunctions on particular variables more prone to yield irregular cuts. 
Nevertheless, 
our results indicate that one could use information about cuts from disjunctions where $|K|=2$ 
to augment and combine those sets of indices yielding irregular cuts 
in the hope of obtaining irregular cuts on disjunctions with more terms.

Fifth, 
the frequency of irregular cuts depends on the structure of the problem, 
which can be higher in larger instances  
but lower when more inequalities are added. 
This comes from observing the results for the many families of instances in the experiments. 
For families with similar size, 
we observe that some have little irregularity in the results (\texttt{mas74} and \texttt{mas76}, all \change{ instances prefixed with \texttt{misc} except \texttt{misc05}, instances prefixed with \texttt{bell}, and instances prefixed with \texttt{fixnet}}), 
some have a moderate level (\texttt{markshare1} and \texttt{markshare2}, \texttt{mod008} and \texttt{mod013}, and \texttt{aflow30a} and \texttt{aflow40b}), 
and some are highly irregular (\change{ instances prefixed with \texttt{set1}  and instances prefixed with \texttt{timtab}}). 
Notably, 
instances in the latter sets are often larger. 
We also observe instances yielding more irregular cuts according to their size 
\change{ among instances prefixed with \texttt{p}} (\texttt{p0033} to \texttt{p0291}), \texttt{stein9} to \texttt{stein27}, and \texttt{vpm1} and \texttt{vpm2}. 
Curiously, however, 
we observe a reduction in the number of \change{strictly} irregular cuts when 
\change{ comparing pairs of instances that represent a same problem. 
For example, instances  \texttt{pp08a}, \texttt{gesa2\_o}, \texttt{gesa3\_o}, and \texttt{qnet1\_o} 
are respectively equivalent to instances \texttt{pp08aCUTS}, \texttt{gesa2}, \texttt{gesa3}, and \texttt{qnet1}. 
However, the latter four have additional constraints strengthening the formulation. 
Among those, the number of strictly irregular cuts only increases from \texttt{qnet1\_o} to \texttt{qnet1},} 
but in that case the number of fractional variables is substantially smaller and that may have affected the results. 
\change{This finding resonates with the example in Figure~\ref{fig:illustration}d, which shows that additional inequalities can make more of the valid cuts regular.} 
We could hypothesize that these additional constraints extend the reach of regular cuts, 
hence making irregular cuts less relevant than before. 
Overall, these results indicate that cuts from disjunctions with more terms would yield more irregular cuts at the root node 
and that they are more likely to do so in larger problems. 

\change{
Finally, irregular cuts are as strong as regular cuts, and thus could complement them if necessary. When comparing the average gap closed or the average Euclidean distance to the fractional solution separated, there are about the same number of cases favoring either type of cut. Figure~\ref{fig:total_gap} shows that, when both regular and irregular cuts are generated for a given problem, part of the optimality gap closed in many of the instances is due exclusively to the irregular cuts. In fact, some irregular cuts may be obtained as regular cuts of higher rank. For example, the irregular cut in Figure~\ref{fig:illustration}c can be obtained as a regular cut of the tightened LP relaxation in Figure~\ref{fig:illustration}b, but the use of previously generated cuts in the LP tableau to generate further intersection cuts may propagate numerical errors from those previous cuts as well. Hence, directly generating those cuts may be preferable in certain cases. 
From Figure~\ref{fig:average}, 
we can 
note that the results for the average Euclidean distance fluctuate at a much smaller scale around the identity line, possibly indicating that outliers to either side could be a starting point to refine the analysis in future work.}

\change{Most of the related computational work has focused on regular cuts, including \citet{MRDP}, which is understandable since those are the cuts that we know how to generate more easily.  
Even among those, \citet{Louveaux} note that general-purpose cut generation is not yet competitive with specific methods, and that using too many rows of the tableau of the LP relaxation is not necessarily better. 
In fact, most work is based on methods that are tailored for the different types of maximal lattice-free sets used to generate intersection cuts~\citep{MultiRowsExp,BasuExp,MRDP,DeyExp}. 
Future work on irregular cuts should probably focus on tailored methods.}

\change{One exception to the focus on regular cuts is the work on cross-cuts by~\citet{CrossExps}, 
which also formulates a CGLP based on non-split disjunctions. 
Interestingly, some ideas in this paper resonate with our observations and findings. 
One of their best approaches is denoted as Cross.def, 
which generates cross cuts that could only be otherwise obtained from split disjunctions as cuts of rank 2. 
Similarly, 
we have argued that irregular cuts can be regular cuts of higher rank, 
even for the same type of disjunction. 
In another approach described in the paper, 
one of the disjunctions is fixed and a second disjunction is chosen to generate a cut through an MILP formulation. Our finding about how simple 3-branch and 4-branch split disjunctions yielding strictly irregular cuts are associated with those yielding strictly irregular cuts among simple 2-branch split disjunctions paves the way to judiciously augment disjunctions when aiming to generate strictly irregular cuts. 
}

\section{Conclusion}

In this paper we have used a mixed-integer formulation to determine when the equivalence between lift-and-project cuts from arbitrary disjunctions and intersection cuts does not hold.  
This method is conveniently used to evaluate the extent to which 
unstrengthened 
$t$-branch split cuts differ from multi-row cuts, 
two types of cuts that have been intensely studied for the past decade. 
When there is no equivalence, the cut is said to be irregular and can only be obtained from irregular CGLP solutions. 

\change{We have focused our analysis on strictly irregular cuts, which we have shown to conveniently exclude irregular cuts that are implied by the set of regular cuts.}  
On the one hand, 
we have found 
that the incidence of irregular cuts varies across different families of instances 
and that many irregular CGLP solutions nevertheless correspond to regular cuts. 
On the other hand, 
the incidence of irregular cuts often 
increases as the problems get larger, the linear relaxation is weaker, and the disjunction has more terms. 
\change{We also found that proportionally less basic regular CGLP solutions characterize split cuts in disjunctions with more terms. 
Interestingly, we observed in our experiments that irregular cuts are on average as strong as regular cuts and that they are often  responsible for part of the optimality gap closed when both types of cuts are generated.  
Hence, irregular cuts could complement their regular counterparts, 
in particular to avoid the numerical instability of higher rank cuts that could otherwise be necessary to solve some problems. In such case, the only irregular cuts that would be relevant are precisely the strictly irregular ones that we analyzed in this study.} 
Furthermore, 
we have observed that 
$3$-branch and $4$-branch split disjunctions yield irregular cuts more often 
when they augment multiple $2$-branch split disjunctions that also yield irregular cuts, 
a result that could be used to judiciously choose disjunctions with more terms from which to generate additional 
cuts.    
\change{Ultimately, given that we would be mostly interested in irregular cuts when generating lift-and-project cuts from arbitrary disjunctions, knowing which disjunctions are more likely to generate irregular cuts can be extremely helpful.}

\paragraph{Acknowledgments} The research of Egon Balas was supported by NSF Grant CMMI-1560828 and ONR contract N00014-18-1-212.

\bibliographystyle{plainnat}
\bibliography{I-CGLP}

\end{document}